\documentclass{article}
\usepackage{arxiv}
\usepackage{amsmath}
\usepackage[utf8]{inputenc} % allow utf-8 input
\usepackage[T1]{fontenc}    % use 8-bit T1 fonts
\usepackage{hyperref}       % hyperlinks
\usepackage{url}            % simple URL typesetting
\usepackage{booktabs}       % professional-quality tables
\usepackage{amsfonts}       % blackboard math symbols
\usepackage{nicefrac}       % compact symbols for 1/2, etc.
\usepackage{microtype}      % microtypography
\usepackage{graphicx}
\usepackage{float}%newly added
\usepackage{afterpage}%newly added
\usepackage{booktabs, calc}%newly added
\usepackage{multicol}
\usepackage{multirow}%newly added
\usepackage{indentfirst}
\usepackage{algorithm}%newly added
\usepackage{algpseudocode}
\usepackage{subcaption}
\usepackage{algorithm} %newly added

\DeclareMathAlphabet\mathbfcal{OMS}{cmsy}{b}{n} % new one i have addeded here
\DeclareMathOperator{\sech}{sech} % new one i have added here
\usepackage{lineno}

\newtheorem{theorem}{Theorem}
\newtheorem{lemma}{Lemma}

\newtheorem{definition}[theorem]{Definition}

\newtheorem{remark}{Remark}

\renewcommand{\Re}{\mathrm{Re} \;}

\newcommand\norm[1]{\left\lVert#1\right\rVert}

\title{Fast and accurate high-order method for high dimensional space-fractional reaction-diffusion equation with general boundary conditions}

\author{
Almushaira Mustafa  \\
School of Mathematics and Statistics\\
Huazhong University of Science and Technology\\
Wuhan 430074, China\\
Department of Mathematics\\
Sana'a University\\
Sana'a, Yemen\\
  \texttt{mstf1985@hust.edu.cn} \\
   \And
Harish Bhatt \\
  Department of Mathematics\\
  Center for Student Success \& Retention\\
  Savannah State University\\
 Savannah, GA 31404 \\
  \texttt{bhatth@savannahstate.edu} \\
}

\begin{document}
\maketitle
\begin{abstract}
%%% Text of abstract
\indent To achieve efficient and accurate long-time integration, we propose a fast, accurate, and stable high-order numerical method for solving fractional-in-space reaction-diffusion equations. The proposed method is explicit in nature and utilize the fourth-order compact finite difference scheme and matrix transfer technique (MTT) in space with FFT-based implementation. A time integration is done through the accurate fourth-order modified exponential time differencing Runge-Kutta scheme. The linear stability analysis and various numerical experiments including two-dimensional (2D) Fitzhugh-Nagumo, Gierer-Meinhardt, Gray-Scott and three-dimensional (3D) Schnakenberg models are presented to demonstrate the accuracy, efficiency and stability of the proposed method.
\end{abstract}

% keywords can be removed

\section{Introduction}
\label{S1}
In this paper, a space-fractional reaction-diffusion equation (SFRDE) of the following form is considered:
\begin{align}\label{Eq1}
\begin{cases}
\dfrac{\partial  u(\mathbf{x},t)}{\partial  t} &= - \kappa (-\Delta)^{\alpha/2} u(\mathbf{x},t) + f(u(\mathbf{x},t),t), \quad  (\mathbf{x},t) \in \Omega\times (0,T], \hspace{0.1cm} 1 < \alpha \leq 2,\\
u(\mathbf{x},0) & = g(\mathbf{x}), \quad \mathbf{x} \in \Omega ,
\end{cases}
\end{align}
where $ \Omega $ is an open bounded domain in $ \mathbb{R}^d $, ($ d=1,2,3 $),  $ \kappa $ represents the diffusion coefficient, $ (-\Delta)^{\alpha/2}  $ is the fractional Laplacian of order $ \alpha $
\cite{ilic2005numerical,yang2010numerical,turner2010use,constantin2016remarks,song2017computing} defined through the eigenfunction expansion on a finite domain (see Definition \ref{D1}) and
$ f(u)  $ is a nonlinear reaction term. The SFRDE \eqref{Eq1} is provided with an initial condition $ g(\mathbf{x}) $ and general boundary conditions such as periodic, homogeneous Dirichlet or Neumann boundary conditions.

Most numerical schemes for solving space-fractional diffusion problems involve applying the finite difference, finite volume, finite element, spectral method and so on for discretizing the fractional operator. Here, we focus our attention on matrix transfer technique (MTT) introduced by Ilic et al. in \cite{ilic2005numerical,ilic2006numerical} to solve space-fractional diffusion equations. Ding and Zhang \cite{2012new} presented fourth-order methods based on the MTT for solving Riesz space-fractional diffusion and advection-dispersion equations. Yang et al. \cite{yang2010numerical} utilized the MTT to solve fractional partial differential equations (FPDEs) with Riesz space-fractional derivative (RSFD). A rational approximation to dense matrix obtained from the MTT was proposed by Aceto and Novati in \cite{aceto2017rational} and its integral representation is approximated with the Guass Jacobi quadrature rule.

There are many existing definitions of the fractional Laplacian operator, for more details, the readers are referred to \cite{pozrikidis2016fractional,duo2017comparative}. It is important to choose an appropriate numerical method of approximation accor to which definition is intended. The basic idea of the MTT is to approximate the fractional Laplacian $ (-\Delta)^{\alpha/2}  $
by the matrix representation $ A^{\alpha/2} $, where $ A $ is a symmetric positive definite matrix obtained from the discretization of the standard Laplace operator subject to the given boundary conditions. The main advantage of this approach is that, it gives a full diagonal representation of the fractional operator, being able to efficiently implement regardless of the fractional power in the problem. An additional advantage is that the application to two and three spatial dimensions is essentially the same as the one dimensional problem.

In this work, we propose a fast, accurate and stable compact exponential time differencing Runge-Kutta method for solving the multi-dimensional SFRDE \eqref{Eq1}. We use a fourth-order compact finite difference method and MTT for the spatial discretization, yiel a diagonalization matrix system of ordinary differential equations (ODEs) whose solution can be explicitly expressed in term of the time integrator. The proposed method can deal with various boundary conditions and take advantage of FFT-based computations. All these techniques are coupled together to produce fast, accurate and stable numerical method.

The rest of the paper is organized as follows. In Section \ref{S2}, we recall an important definition and lemma. In Section \ref{S3}, we discuss 1D, 2D and 3D versions of the proposed spatial discretization method with separate discussions for different types of boundary conditions. The time discretization is introduced in Section \ref{S4} and its linear stability analysis is presented. Furthermore, the 3D algorithm of the proposed method is provided with periodic boundary conditions. Extensive numerical experiments with applications are reported in Section \ref{S5} in order to numerically demonstrate accuracy, efficiency and stability of the proposed method. Some conclusion remarks are made in Section \ref{S6}.

\section{Preliminaries}
\label{S2}
In this section, the following definition and lemma which will be useful throughout this paper are recalled.

\begin{definition}\label{D1}
	For $ n,m,l=0,1,... $,
	 let the Laplacian $  (-\Delta)$ has a complete set of orthogonal eigenfunction $ \psi_n $, $ \psi_{n,m} $, or $ \psi_{n,m,l} $ corresponding to eigenvalues $ \lambda_n $, $ \lambda_{n,m} $, or $ \lambda_{n,m,l} $, respectively, on a bounded region $ \Omega $ with periodic, homogeneous Dirichlet or homogeneous Neumann boundary conditions, then $  (-\Delta)^{\alpha/2}  $ is defined by
	 \begin{align*}
	 \begin{split}
	  (-\Delta)^{\alpha/2} u & = \sum_{n=0}^{\infty} \hat{u}_n \lambda_n^{\alpha/2} \psi_n, \hspace{2.3cm} d=1,\\
	   (-\Delta)^{\alpha/2} u & =\sum_{n=0}^{\infty} \sum_{m=0}^{\infty} \hat{u}_{n,m} \lambda_{n,m}^{\alpha/2} \psi_{n,m}, \hspace{1.3cm} d=2,\\
	   (-\Delta)^{\alpha/2} u & =\sum_{n=0}^{\infty} \sum_{m=0}^{\infty} \sum_{l=0}^{\infty} \hat{u}_{n,m,l} \lambda_{n,m,l}^{\alpha/2} \psi_{n,m,l}, \hspace{.5cm} d=3,
	 \end{split}
	 \end{align*}
	 in which $ u $ can be expressed as follows
	 \begin{align*}
	 \begin{split}
	 u & = \sum_{n=0}^{\infty} \hat{u}_n \psi_n, \quad \text{such that}\quad \sum_{n=0}^{\infty} |\hat{u}_n|^2 |\lambda_n|^\alpha < \infty,  \hspace{3.6cm} d=1,\\
	 u & = \sum_{n=0}^{\infty} \sum_{m=0}^{\infty} \hat{u}_{n,m} \psi_{n,m}, \quad \text{such that}\quad \sum_{n=0}^{\infty} \sum_{m=0}^{\infty} |\hat{u}_{n,m}|^2 |\lambda_{n,m}|^\alpha < \infty,  \hspace{1.7cm} d=2,\\
	 u & = \sum_{n=0}^{\infty} \sum_{m=0}^{\infty} \sum_{l=0}^{\infty} \hat{u}_{n,m,l} \psi_{n,m,l}, \quad \text{such that} \quad\sum_{n=0}^{\infty} \sum_{m=0}^{\infty} \sum_{l=0}^{\infty}  |\hat{u}_{n,m,l}|^2 |\lambda_{n,m,l}|^\alpha < \infty,  \hspace{0.2cm} d=3.
	  \end{split}
	 \end{align*}
\end{definition}

\begin{lemma}[\cite{zhuang2009numerical}]
Suppose $ A $ be a positive definite matrix, and $ A=P \Lambda P^{-1} $, where $ P $ is an orthogonal matrix and $ \Lambda $ is a diagonal matrix with diagonal entries being the eigenvalues of $ A $, then for arbitrary real $ \alpha $, $ A^{\alpha} = P \Lambda^{\alpha} P^{-1}$ may be uniquely determined by $ A $ and $ \alpha $.
\end{lemma}

\section{Matrix transfer technique}
\label{S3}
In this section, a fourth-order compact difference scheme for discretizing the fractional operator $  (-\Delta)^{\alpha/2} $ is introduced and derived a fast algorithm via FFT-based implementation.

By using MTT developed by Ding and Zhang in \cite{2012new} for discretizing the fractional operator $  (-\Delta)^{\alpha/2} $ with a uniform mesh of step size $ h $ in each spatial direction we obtain
\begin{equation}\label{MTT}
(-\Delta)^{\alpha/2} u \approx  \mathcal{T}^{\alpha/2} u,
\end{equation}
where $ \mathcal{T}^{\alpha/2} $ is constructed from the eigenvalues and eigenvectors of the matrix representation of the standard Laplacian $  \mathcal{T} $ based on given boundary conditions. That is, the matrix $   \mathcal{T}^{\alpha/2} $ will be represented in the following subsections need not be formed explicitly.

Here, we use a fourth-order compact difference scheme for the Laplace operator %$  (-\Delta) $
\begin{equation}\label{Lap}
 -\Delta \approx -h^{-2} (1 + \frac{1}{12} \delta^2_x)^{-1} \delta^2_x + \mathcal{O}(h^4)
\end{equation}
where $ \delta^2_x $ is the second-order central difference operator.
\subsection{One-dimensional case }
Let us consider the model \eqref{Eq1} in an open domain in the one-dimensional space $ \Omega=\{x_a < x <x_b\} $. We denote $ \{x_n\}_{n=0}^{N} $ as the uniform grid points and $ h=\frac{x_b-x_a}{N} $ as the mesh size. Let $ u_n=u_n(t)\approx u_n(x_n,t) $ for $ 0\leq n \leq N $ denotes the numerical solution. Denote the solution vector as $ \mathbf{u}=(u_0,\cdots,u_N)^T $.

\subsubsection{The case of periodic boundary condition}
For the periodic boundary condition,
the matrix representation of standard Laplacian based on a fourth-order compact difference discretization reads:
\begin{equation*}
-\Delta u(x,t) =  \left( A^{-1}B \right) \mathbf{u}
\end{equation*}
where $ A_{N\times N}$, $ B_{N\times N} $ are the tridiagonal matrices as follows
\begin{equation*}
A = 
\begin{bmatrix}
\frac 56 &\frac {1}{12}  & 0 &\cdots  &0  &\frac {1}{12}  \\
\frac {1}{12}&\frac{5}{6}  &\frac{1}{12} &0 &\cdots  &0 \\
0&\frac{1}{12}  & \frac 56   &\frac{1}{12}  &0  &\cdots\\
\vdots&\ddots  &\ddots  &\ddots   &  &  \vdots \\ 
\cdots&0 & \frac{1}{12}  & \frac 56   &\frac{1}{12}&0   \\ 
0&  \cdots&0 & \frac{1}{12}  & \frac 56   &\frac{1}{12}   \\
\frac {1}{12}&0&  \cdots&0 & \frac{1}{12}  & \frac 56 
\end{bmatrix}
\end{equation*}

\begin{equation*}
B= h^{-2}
\begin{bmatrix}
2 &-1 & 0 &\cdots  &0  &-1 \\
-1&2  &-1&0 &\cdots  &0 \\
0&-1 & 2  &-1 &0  &\cdots\\
\vdots&\ddots  &\ddots  &\ddots   &  &  \vdots \\ 
\cdots&0 & -1 & 2  &-1&0   \\ 
0&  \cdots&0 & -1  & 2   &-1  \\
-1&0&  \cdots&0 & -1  & 2
\end{bmatrix}
\end{equation*}
and 
\[
\mathbf{u}=(u_1,\cdots,u_N)^T .
\]
Now, if $  \mathcal{T} = A^{-1}B $ is a matrix representation of a differential operator $ -\Delta $, then the matrix representation of the fractional operator $ (-\Delta )^{\alpha/2} $ can be given by $   \mathcal{T}^{\alpha/2}  $.
%where $ \Lambda $ is the diagonal eigenvalues matrix and $ P $ is the eigenvectors matrix.
% h^{-\alpha} (P \Lambda^{\alpha/2} P^{-1})
Since the matrices $ A $ and $ B $ are symmetric positive definite, they can be diagonalized as \cite{ju2015fast,zhu2016fast,huang2019afast,huang2019fast,alzahrani2019fourth}
\begin{align}\label{MP}
\mathcal{T} = A^{-1}B = P \Lambda P^{-1}= \left(  P \Lambda P^{-1} \right)^{T} = P^{-1} \Lambda P,
\end{align}
where $ \Lambda $ is the diagonal eigenvalues matrix, that is,
\begin{equation}
\Lambda =diag(\lambda_1,\cdots,\lambda_N), \quad 
 \lambda_n=\frac{4 \sin^2(\frac{(n-1) \pi }{2 N}) }{h^2\left( 1-\frac 13 \sin^2(\frac{(n-1) \pi }{2 N}) \right) }, \quad
 1\leq n \leq N,
\end{equation}
and $ P_{N\times N}$ is the eigenvectors matrix, that is,
\begin{equation}\label{PP}
P_{i,j}=\exp\left( \frac{-\sqrt{-1} (i-1)(j-1) 2 \pi }{N}\right) , \quad 
i,j=1,\cdots,N.
\end{equation}
Using the MTT, the matrix representation of fractional Laplacian reads 
\begin{equation}\label{flp}
(-\Delta)^{\alpha/2} u(x,t) =\left(  P^{-1} \Lambda^{\alpha/2} P \right) \mathbf{u}.
\end{equation}
Substituting \eqref{flp} into \eqref{Eq1} yields the semi-discretization in space of the one-dimensional SFRDE:
\begin{align}\label{S7}
\dfrac{\partial  \mathbf{u}}{\partial  t} &= - \kappa  \left(  P^{-1}  \Lambda^{\alpha/2} P \right) \mathbf{u}  + \mathbf{f}( \mathbf{u},t),
\end{align}
where $ \mathbf{f}(\mathbf{u},t)=\left( f(u_1,t),\cdots,  f(u_N,t) \right)^T  $.
If $ N $ is large, the computational cost of the matrix-vector product in \eqref{flp} is extremely high. Hence, we will propose a fast algorithm to efficiently compute matrix-vector product in \eqref{flp}. 
From the definition of $ P $ in \eqref{PP}, it is pointed out in \cite{ju2015fast} that the product $ P \mathbf{u} $ is equivalent to the Discrete Fourier transform (DFT) \cite{van1992computational} of the vector $ \mathbf{u} $, namely, $ P \mathbf{u}=\mathcal{F}(\mathbf{u}) $, ( $\mathcal{F}  $ stands for DFT ), while $ P^{-1} \mathbf{u} $ can be obtained via the inverse DFT of $ \mathbf{u}  $, that is, $ P^{-1} \mathbf{u}  = \mathcal{F}^{-1} (\mathbf{u} )$. 

Multiplying $ P $ from the left hand side in \eqref{S7}, we immediately obtain
%Therefore, we can rewrite the semi-discretization scheme in \eqref{S7} as follows
\begin{align}\label{S8}
\dfrac{\partial  \mathcal{F}( \mathbf{u})}{\partial  t} &= - \kappa  \Lambda^{\alpha/2}  \mathcal{F}(\mathbf{u})  + \mathcal{F} (\mathbf{f}( \mathbf{u},t)) .
\end{align}
Thus, the computational complexity of evaluating the matrix-vector product in \eqref{flp} can be reduced from $ \mathcal{O}(N^2) $ to $ \mathcal{O}(N \log(N)) $.

\subsubsection{The case of homogenous Dirichlet boundary condition}

For the homogenous Dirichlet boundary condition, the matrix representation of standard Laplacian reads as:
\begin{equation*}
-\Delta u(x,t) = \left( A^{-1}B \right) \mathbf{u},
\end{equation*}
where $ A_{(N-1)\times (N-1)}$, $ B_{(N-1)\times (N-1)} $ are the tridiagonal matrices as follows
\begin{equation*}
A = 
\begin{bmatrix}
\frac 56 &\frac {1}{12}  & 0 &\cdots  &0  &0 \\
\frac {1}{12}&\frac{5}{6}  &\frac{1}{12} &0 &\cdots  &0 \\
0&\frac{1}{12}  & \frac 56   &\frac{1}{12}  &0  &\cdots\\
\vdots&\ddots  &\ddots  &\ddots   &  &  \vdots \\ 
\cdots&0 & \frac{1}{12}  & \frac 56   &\frac{1}{12}&0   \\ 
0&  \cdots&0 & \frac{1}{12}  & \frac 56   &\frac{1}{12}   \\
0&0&  \cdots&0 & \frac{1}{12}  & \frac 56 
\end{bmatrix}
\end{equation*}

\begin{equation*}
B= h^{-2}
\begin{bmatrix}
2 &-1 & 0 &\cdots  &0  &0\\
-1&2  &-1&0 &\cdots  &0 \\
0&-1 & 2  &-1 &0  &\cdots\\
\vdots&\ddots  &\ddots  &\ddots   &  &  \vdots \\ 
\cdots&0 & -1 & 2  &-1&0   \\ 
0&  \cdots&0 & -1  & 2   &-1  \\
0&0&  \cdots&0 & -1  & 2
\end{bmatrix}
\end{equation*}
and 
\[
\mathbf{u}=(u_1,\cdots,u_{N-1})^T .
\]

Similarly, the matrix representation of fractional Laplacian reads: 
\begin{equation*}
(-\Delta)^{\alpha/2} u(x,t) =\left(  P^{-1} \Lambda^{\alpha/2} P \right) \mathbf{u},
\end{equation*}
where 
\begin{equation}
\Lambda =diag(\lambda_1,\cdots,\lambda_{N-1}), \quad 
\lambda_n=\frac{4 \sin^2(\frac{n \pi }{2 N}) }{h^2\left( 1-\frac 13 \sin^2(\frac{n \pi }{2 N})\right)  }, \quad
1\leq n \leq N-1,
\end{equation}
and
\begin{equation}\label{PD}
P_{i,j}=\sin\left( \frac{i j \pi }{N}\right) , \quad 
i,j=1,\cdots,N-1.
\end{equation}
From the definition of $ P $ in \eqref{PD}, 
it is noted that the operations $ P \mathbf{u} $ and $ P^{-1} \mathbf{u} $ are exactly the Discrete Sine Transform (DST) and the inverse DST respectively, and can be efficiently evaluated by FFT-based algorithm \cite{van1992computational}.

\subsubsection{The case of homogenous Neumann boundary condition}

For the homogenous Neumann boundary condition, the matrix representation of standard Laplacian reads: 
\begin{equation*}
-\Delta u(x,t) =\left( A^{-1}B \right) \mathbf{u},
\end{equation*}
where $ A_{(N+1)\times (N+1)}$, $ B_{(N+1)\times (N+1)} $ are the tridiagonal matrices as follows
\begin{equation*}
A = 
\begin{bmatrix}
\frac 56 &\frac {1}{6}  & 0 &\cdots  &0  &0 \\
\frac {1}{12}&\frac{5}{6}  &\frac{1}{12} &0 &\cdots  &0 \\
0&\frac{1}{12}  & \frac 56   &\frac{1}{12}  &0  &\cdots\\
\vdots&\ddots  &\ddots  &\ddots   &  &  \vdots \\ 
\cdots&0 & \frac{1}{12}  & \frac 56   &\frac{1}{12}&0   \\ 
0&  \cdots&0 & \frac{1}{12}  & \frac 56   &\frac{1}{12}   \\
0&0&  \cdots&0 & \frac{1}{6}  & \frac 56 
\end{bmatrix}
\end{equation*}

\begin{equation*}
B=  h^{-2}
\begin{bmatrix}
2 &-2 & 0 &\cdots  &0  &0\\
-1&2  &-1&0 &\cdots  &0 \\
0&-1 & 2  &-1 &0  &\cdots\\
\vdots&\ddots  &\ddots  &\ddots   &  &  \vdots \\ 
\cdots&0 & -1 & 2  &-1&0   \\ 
0&  \cdots&0 & -1  & 2   &-1  \\
0&0&  \cdots&0 & -2  & 2
\end{bmatrix}
\end{equation*}
and 
\[
\mathbf{u}=(u_0,\cdots,u_{N})^T .
\]

Similarly, the matrix representation of fractional Laplacian reads: 
\begin{equation*}
(-\Delta)^{\alpha/2} u(x,t) =\left(  P^{-1}  \Lambda^{\alpha/2} P\right) \mathbf{u},
\end{equation*}
where 
\begin{equation}
\Lambda =diag(\lambda_1,\cdots,\lambda_{N-1}), \quad 
\lambda_n=\frac{4 \sin^2(\frac{n \pi }{2 N}) }{h^2\left( 1-\frac 13 \sin^2(\frac{n \pi }{2 N}) \right) }, \quad
0\leq n \leq N,
\end{equation}
and
\begin{equation}\label{PN}
P_{i,j}=\cos\left( \frac{i j \pi }{N}\right) , \quad 
i,j=0,\cdots,N.
\end{equation}
In order to reduce the computational complexity of $ P \mathbf{u}$ or $ P^{-1} \mathbf{u}$ from $ \mathcal{O}(N^2) $ to $ \mathcal{O}(N \log(N)) $ by a FFT-based fast algorithm, we follow the idea in \cite{ju2015fast}. To this end, let $ \mathbf{u} =(u_0,\cdots,u_N)^T $ be any vector of size $ N+1 $ and similarly define its reflection vector of size $ N-1 $ as $ \mathbf{v} =(u_{N-1},\cdots,u_1)^T $, and apply a DFT to $ (\mathbf{u} , \mathbf{v})^T $ and then take the first $ N+1 $ components of the result which is equivalent to operation $ P  \mathbf{u}$ and also apply an inverse DFT to $ (\mathbf{u} , \mathbf{v})^T $ and then take the first $ N+1 $ components of the result which is equivalent to operation $ P^{-1}  \mathbf{u}$.

\subsection{Two-dimensional case }

For $ d=2 $, we suppose $ \Omega=\{x_a < x < x_b, y_a < y < y_b\} $. With loss of generality, we partition the spatial domain $ \Omega $ by a square grid which is uniform in each direction, i.e., 
$ h=\frac{x_b-x_a}{N} =\frac{y_b-y_a}{N} $. Mesh points are defined as $ (x_n,y_m)= (x_a + n h, y_a +m h) $ for $ n,m=0,\cdots,N $. Let $ u_{n,m} = u_{n,m}(t)  $ represent the numerical approximation of the solution $ u(x_n,y_m,t) $, for $ n,m=0,\cdots,N $. Denote the 2D solution array as $ U=\{u_{n,m}\}_{(N+1,N+1)} $.

%\subsubsection{The case of periodic boundary condition}

For the periodic boundary condition, following the same arguments of the one-dimensional case, the discretization of fractional Laplacian can be expressed as:
\begin{equation}\label{EP}
(-\Delta)^{\alpha/2} u(x,y,t) = \left(   \mathcal{T}^{\alpha/2} _x U + U (\mathcal{T}^{\alpha/2} _y)^T   \right) ,
\end{equation}
where the matrices $  \mathcal{T}_x = A^{-1} B $ and $  \mathcal{T}_y = A^{-1} B $ are given in \eqref{MP}.
% $ I \in \mathbb{R}^{N\times N}$ represents the identity matrix and $ \otimes $ represents the Kronecker product of two matrices. 
 %We can rewrite \eqref{EP} in term of the solution array $ U $. To this end,  
 Define the special operators $ \raisebox{0.8pt}{\textcircled{\raisebox{0.5pt} {x}}}  $ and
  $  \raisebox{0.8pt}{\textcircled{\raisebox{0.5pt} {y}}} $ as follows
 %For any $ N\times N $ matrix $ U $, we define 
 $  \mathcal{T}_x \raisebox{0.8pt}{\textcircled{\raisebox{0.5pt} {x}}} U = \mathcal{T}_x U$ and  $  \mathcal{T}_y \raisebox{0.8pt}{\textcircled{\raisebox{0.5pt} {y}}} U = U \mathcal{T}^T_y $, that is
 \begin{align*}
 \left( \mathcal{T}_x \raisebox{0.8pt}{\textcircled{\raisebox{0.5pt} {x}}} U\right)_{n,m} = \sum_{k=1}^{N}
 (\mathcal{T}_x)_{n,k} u_{k,m} ,\quad 
  \left( \mathcal{T}_y \raisebox{0.8pt}{\textcircled{\raisebox{0.5pt} {y}}} U\right)_{n,m} = \sum_{k=1}^{N}
 (\mathcal{T}_y)_{n,k} u_{n,k} .
 \end{align*}
Note that these two operators are commutative, i.e., $ \mathcal{T}_x \raisebox{0.8pt}{\textcircled{\raisebox{0.5pt} {x}}} \mathcal{T}_y \raisebox{0.8pt}{\textcircled{\raisebox{0.5pt} {y}}} U = \mathcal{T}_y \raisebox{0.8pt}{\textcircled{\raisebox{0.5pt} {y}}} \mathcal{T}_x \raisebox{0.8pt}{\textcircled{\raisebox{0.5pt} {x}}}  U$.
Then we can rewrite \eqref{EP} in the following compact representation
\begin{equation}
(-\Delta)^{\alpha/2} u(x,y,t) =h^{-\alpha} \left(   \mathcal{T}^{\alpha/2} _x \raisebox{0.8pt}{\textcircled{\raisebox{0.5pt} {x}}} U + \mathcal{T}^{\alpha/2} _y \raisebox{0.8pt}{\textcircled{\raisebox{0.5pt} {y}}} U  \right).
\end{equation}
Define another operator $ \odot $ for element by element multiplication two arrays of same sizes as
\[
(A \odot B)_{i,j} = (B \odot A)_{i,j} = A_{i,j} B_{i,j}.
\]
Then, the fractional Laplacian can be reformulated as
\begin{equation}\label{E16}
(-\Delta)^{\alpha/2} u(x,y,t) =  P^{-1}_y \raisebox{0.8pt}{\textcircled{\raisebox{0.5pt} {y}}} P^{-1}_x \raisebox{0.8pt}{\textcircled{\raisebox{0.5pt} {x}}} \left( \Lambda^{\alpha/2} \odot \left( 
P_y \raisebox{0.8pt}{\textcircled{\raisebox{0.5pt} {y}}} P_x \raisebox{0.8pt}{\textcircled{\raisebox{0.5pt} {x}}} U \right) \right),
\end{equation}
where 
$ \Lambda_{n,m}= h^{-2}\left( \frac{4 \sin^2(\frac{(n-1) \pi }{2 N}) }{1-\frac 13 \sin^2(\frac{(n-1) \pi }{2 N}) } +
\frac{4 \sin^2(\frac{(m-1) \pi }{2 N}) }{1-\frac 13 \sin^2(\frac{(m-1) \pi }{2 N}) }\right) $, for $ n,m=1,\cdots,N $,
$ P_x $ and $ P_y $ are orthogonal matrices be given in \eqref{PP} and $ U=\{u_{n,m}\}_{(N,N)} $. 

For the homogenous Dirichlet boundary condition, the fractional Laplacian can be similarly expressed as:
\begin{equation}\label{E17}
(-\Delta)^{\alpha/2} u(x,y,t) =  P^{-1}_y \raisebox{0.8pt}{\textcircled{\raisebox{0.5pt} {y}}} P^{-1}_x \raisebox{0.8pt}{\textcircled{\raisebox{0.5pt} {x}}} \left(  \Lambda^{\alpha/2} \odot \left( 
P_y \raisebox{0.8pt}{\textcircled{\raisebox{0.5pt} {y}}} P_x \raisebox{0.8pt}{\textcircled{\raisebox{0.5pt} {x}}} U \right) \right),
\end{equation}
where 
$ \Lambda_{n,m}=h^{-2} \left(  \frac{4 \sin^2(\frac{n \pi }{2 N}) }{1-\frac 13 \sin^2(\frac{n \pi }{2 N}) } +
\frac{4 \sin^2(\frac{m \pi }{2 N}) }{1-\frac 13 \sin^2(\frac{m \pi }{2 N}) }\right) $, for $ n,m=1,\cdots,N-1 $,
$ P_x $ and $ P_y $ are orthogonal matrices be given in \eqref{PD} and $ U=\{u_{n,m}\}_{(N-1,N-1)} $.

For the homogenous Neumann boundary condition, the fractional Laplacian can be similarly expressed as:
\begin{equation}\label{E18}
(-\Delta)^{\alpha/2} u(x,y,t) =  P^{-1}_y \raisebox{0.8pt}{\textcircled{\raisebox{0.5pt} {y}}} P^{-1}_x \raisebox{0.8pt}{\textcircled{\raisebox{0.5pt} {x}}} \left(\Lambda^{\alpha/2} \odot \left( 
P_y \raisebox{0.8pt}{\textcircled{\raisebox{0.5pt} {y}}} P_x \raisebox{0.8pt}{\textcircled{\raisebox{0.5pt} {x}}} U \right) \right),
\end{equation}
where 
$ \Lambda_{n,m}= h^{-2} \left( \frac{4 \sin^2(\frac{n \pi }{2 N}) }{1-\frac 13 \sin^2(\frac{n \pi }{2 N}) } +
\frac{4 \sin^2(\frac{m \pi }{2 N}) }{1-\frac 13 \sin^2(\frac{m \pi }{2 N}) }\right) $, for $ n,m=0,\cdots,N $,
$ P_x $ and $ P_y $ are orthogonal matrices be given in \eqref{PN} and $ U=\{u_{n,m}\}_{(N+1,N+1)} $.
%%%

Plugging the above equations \eqref{E16}, \eqref{E17} or \eqref{E18} into \eqref{Eq1} and multiplying $  P_x  $ from the left hand side and $  P_y^{T}  $  from the right hand side, we obtain the semi-discretization in space of the two-dimensional SFRDE:
\begin{align}\label{E19}
\dfrac{\partial  \mathcal{F}(U)}{\partial  t} &= - \kappa   \Lambda^{\alpha/2}  \mathcal{F}(U)+ \mathcal{F}(f(U,t)),
\end{align}
where $ \mathcal{F}(U)= P_y \raisebox{0.8pt}{\textcircled{\raisebox{0.5pt} {y}}} P_x \raisebox{0.8pt}{\textcircled{\raisebox{0.5pt} {x}}} U$ and $ \mathcal{F}(f(U,t)) =  P_y \raisebox{0.8pt}{\textcircled{\raisebox{0.5pt} {y}}} P_x \raisebox{0.8pt}{\textcircled{\raisebox{0.5pt} {x}}} f(U,t) $.

We note that the operation $ P_y \raisebox{0.8pt}{\textcircled{\raisebox{0.5pt} {y}}} P_x \raisebox{0.8pt}{\textcircled{\raisebox{0.5pt} {x}}} U $ can be deficiency computed by FFT-based computation, thus the computational complexity can be reduced from $ \mathcal{O}(N^3) $ to $ \mathcal{O}(N^2 \log(N)) $.

\subsection{Three-dimensional case }

%Similar to two-dimensional case, 
For $ d=3 $, let $ \Omega=\{x_a < x< x_b, y_a < y < y_b, z_a < z < z_b\} $. For simplicity, we denote $ h$
as the spatial space and $ N $ as the number of grid intervals in each direction. Set $ u_{n,m,l}= u_{n,m,l}(t) \approx u(x_n,y_m,z_l,t) $ for $ n,m,l=0,\cdots,N $. Denote the unknowns as a three-dimensional array $ U=\{u_{n,m,l}\}_{(N+1,N+1,N+1)} $.
%$  U=(u_{n,m,l}) \in \mathbb{R}^{(N+1)\times (N+1)\times(N+1)}$.
With similar process, we can write the fractional Laplacian as the following compact representation
\begin{equation}
(-\Delta)^{\alpha/2} u(x,y,z,t) = \left(   \mathcal{T}^{\alpha/2} _x \raisebox{0.8pt}{\textcircled{\raisebox{0.5pt} {x}}} U + \mathcal{T}^{\alpha/2} _y \raisebox{0.8pt}{\textcircled{\raisebox{0.5pt} {y}}} U + \mathcal{T}^{\alpha/2} _z \raisebox{0.8pt}{\textcircled{\raisebox{0.5pt} {z}}} U  \right),
\end{equation}
where the spacial operator are redefined as
\begin{equation*}
\left( \mathcal{T}_x \raisebox{0.8pt}{\textcircled{\raisebox{0.5pt} {x}}} U\right)_{n,m,l} = \sum_{k=0}^{N}
(\mathcal{T}_x)_{n,k} u_{k,m,l}, \quad
\left( \mathcal{T}_y \raisebox{0.8pt}{\textcircled{\raisebox{0.5pt} {y}}} U\right)_{n,m,l} = \sum_{k=0}^{N}
(\mathcal{T}_y)_{m,k} u_{n,k,l}, \quad 
\left( \mathcal{T}_z \raisebox{0.8pt}{\textcircled{\raisebox{0.5pt} {z}}} U\right)_{n,m,l} = \sum_{k=0}^{N}
(\mathcal{T}_z)_{l,k} u_{n,m,k}.
\end{equation*}

Following the similar analysis in two dimensional case, the fractional Laplacian can be reformulated as follows
 \begin{equation}
 (-\Delta)^{\alpha/2} u(x,y,z,t) =  P^{-1}_z \raisebox{0.8pt}{\textcircled{\raisebox{0.5pt} {z}}}  P^{-1}_y \raisebox{0.8pt}{\textcircled{\raisebox{0.5pt} {y}}} P^{-1}_x \raisebox{0.8pt}{\textcircled{\raisebox{0.5pt} {x}}} \left( \Lambda^{\alpha/2} \odot \left( 
 P_z \raisebox{0.8pt}{\textcircled{\raisebox{0.5pt} {z}}} P_y \raisebox{0.8pt}{\textcircled{\raisebox{0.5pt} {y}}} P_x \raisebox{0.8pt}{\textcircled{\raisebox{0.5pt} {x}}} U \right) \right),
 \end{equation}
 where 
 \begin{align*}
 \Lambda_{n,m,l} &=h^{-2} \left( \frac{4 \sin^2(\frac{(n-1) \pi }{2 N}) }{1-\frac 13 \sin^2(\frac{(n-1) \pi }{2 N}) } +
 \frac{4 \sin^2(\frac{(m-1) \pi }{2 N}) }{1-\frac 13 \sin^2(\frac{(m-1) \pi }{2 N}) } + 
 \frac{4 \sin^2(\frac{(l-1) \pi }{2 N}) }{1-\frac 13 \sin^2(\frac{(l-1) \pi }{2 N}) } \right) , \quad n,m,l=1,\cdots,N,\\
  & \hspace{1cm}  (\text{in the case of periodic boundary conditions})\\
 \Lambda_{n,m,l} &= h^{-2} \left( \frac{4 \sin^2(\frac{n \pi }{2 N}) }{1-\frac 13 \sin^2(\frac{n \pi }{2 N}) } +
 \frac{4 \sin^2(\frac{m \pi }{2 N}) }{1-\frac 13 \sin^2(\frac{m \pi }{2 N}) } + 
 \frac{4 \sin^2(\frac{l \pi }{2 N}) }{1-\frac 13 \sin^2(\frac{l\pi }{2 N}) }\right) , \quad n,m,l=1,\cdots,N-1,\\
  & \hspace{1cm}  (\text{in the case of homogenous Dirichlet boundary conditions})\\
 \Lambda_{n,m,l} &=h^{-2} \left( \frac{4 \sin^2(\frac{n \pi }{2 N}) }{1-\frac 13 \sin^2(\frac{n \pi }{2 N}) } +
 \frac{4 \sin^2(\frac{m \pi }{2 N}) }{1-\frac 13 \sin^2(\frac{m \pi }{2 N}) } + 
 \frac{4 \sin^2(\frac{l \pi }{2 N}) }{1-\frac 13 \sin^2(\frac{l\pi }{2 N}) }\right) , \quad n,m,l=0,\cdots,N,\\
 & \hspace{1cm}  (\text{in the case of homogenous Neumann boundary conditions})
 \end{align*}
and $ P_x $, $ P_y $ and $ P_z $ are orthogonal matrices be given in \eqref{PP} for periodic boundary conditions, be defined in \eqref{PD} for homogenous Dirichlet boundary conditions or be defined in \eqref{PN} for homogenous Neumann boundary conditions.

The following is the semi-discretization in space of the three-dimensional SFRDE:
\begin{align}\label{E22}
\dfrac{\partial  \mathcal{F}(U)}{\partial  t} &= - \kappa   \Lambda^{\alpha/2}  \mathcal{F}(U)+ \mathcal{F}(f(U,t)),
\end{align}
where $ \mathcal{F}(U)= P_z \raisebox{0.8pt}{\textcircled{\raisebox{0.5pt} {z}}} P_y \raisebox{0.8pt}{\textcircled{\raisebox{0.5pt} {y}}} P_x \raisebox{0.8pt}{\textcircled{\raisebox{0.5pt} {x}}} U$ and $ \mathcal{F}(f(U,t)) =  P_z \raisebox{0.8pt}{\textcircled{\raisebox{0.5pt} {z}}} P_y \raisebox{0.8pt}{\textcircled{\raisebox{0.5pt} {y}}} P_x \raisebox{0.8pt}{\textcircled{\raisebox{0.5pt} {x}}} f(U,t) $.

\begin{remark}
The overall computational cost of the proposed method can be reduced from $ \mathcal{O}(N^{d+1}) $ to $ \mathcal{O}(N^{d} \log(N)) $ for $ d=1,2, $ or $ 3 $, per time step by using FFT-based fast calculation.
\end{remark}

\section{Exponential time differencing procedure}
\label{S4}

  Consider the following nonlinear initial boundary value problem
  \begin{equation}\label{ODE}
  U_t + \mathcal{A} U= f(U,t), \in \Omega, t=(0,T], 
  \end{equation}
  where $ \mathcal{A}  $ represents a spatial discretization of Laplacian operator.
 Let $ \tau $ be the temporal step size, then using a variation of constant formula, then the exact solution of \eqref{ODE} is the following recurrence formula
%  \begin{align}\label{Eq5}
%  \dfrac{\partial \hat{u}(\mathbf{x},t) }{\partial  t} &= - \mathcal{A} \hat{u}(\mathbf{x},t)   +\hat{f}( u,t),
%  \end{align}
  \begin{align}\label{Eq5}
U(t_{k+1})&= e^{- \tau \mathcal{A} } U(t_k )   + \tau \int_{0}^{1} e^{- \tau \mathcal{A} (1-s)} f( U(t_k+s \tau ),t_k+s \tau) ds.
  \end{align}
  %where $ \hat{u}=\mathcal{F}(u) $ and $ \hat{f}(u,t)=\mathcal{F}(f(u,t)) $.
  %The various exponential time differencing ( ETD) methods come from how one approximates the integral in above expression. 
  %The ETD schemes solve the linear part exactly and then explicitly approximates the integral part by polynomial approximation.
  
\begin{figure}[H]
 	\centering
 	\includegraphics[width=1.1\linewidth, height=0.3\textheight]{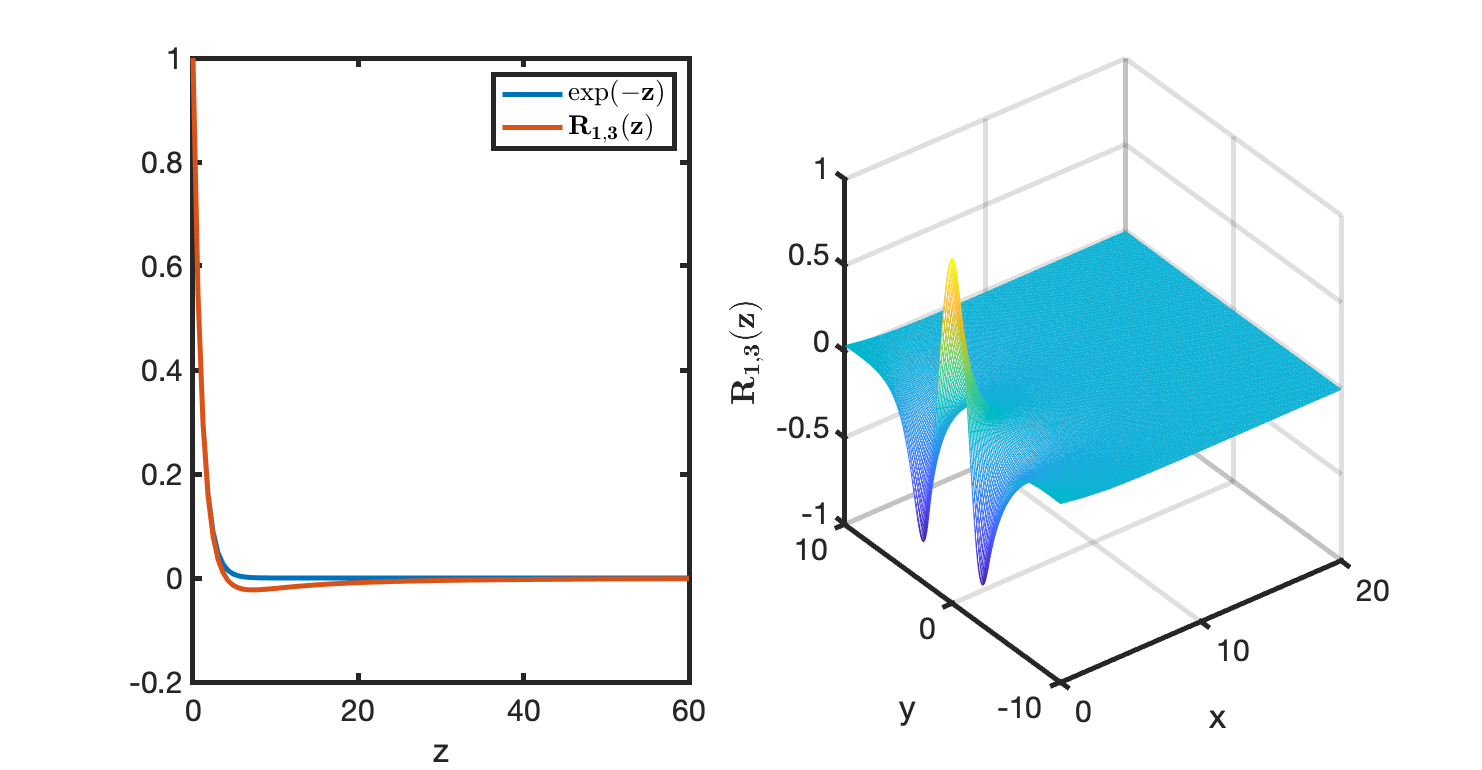}
 	\caption{Left label: Behavior of the functions $ e^{-z} $ and $ R_{1,3}(z) $ for $ z\in [0,60] $.  Right label: Behavior of the function $ R_{1,3}(z) $ approximation of $ e^{-z} $ for $ z=(x,y)\in [0,20]\times [-10,10] $.}
 	\label{fig:r13}
 \end{figure}

  The various exponential time differencing (ETD) methods come into picture from how one approximates the integral term. Here, we utilize a fourth-order time stepping scheme as considered by \cite{bhatt2016compact} which is based on fourth-order Pad\'{e} (1, 3) approximation to $ e^{-z} $. 
 % We use the notation $ R_{1,3}(z) $ for Pad\'{e} (1, 3) approximation to $ e^{-z} $ and $ \tilde{R}_{1,3}(z) $ for Pad\'{e} (1, 3) approximation to $ e^{-z/2} $.
  The Pad\'{e} approximations are known rational approximations. In particular, Pad\'{e} (1, 3) approximation to $ e^{-z} $ is given by
  \begin{align}
  R_{1,3}(z) = \frac{24- 6 z}{24 + 18 z + 6 z^2 +z^3},
  \end{align}
 where the notation $ R_{1,3}(z) $ is used for Pad\'{e} (1, 3) approximation to $ e^{-z} $.
 Fig. \ref{fig:r13} (left) shows the behavior of the exponential function $ e^{-z} $ and $ R_{1,3}(z) $ for $ z\in [0,60]$, while Fig. \ref{fig:r13} (right) shows the behavior of the function $ R_{1,3}(z) $ for $ z=(x,y)\in [0,20]\times [-10,10] $.
 
 Following the process in \cite{bhatt2016compact}, the ETDRK4-P13 method reads as:
 \begin{align}\label{ETDRK4P13}
 \begin{split}
 U^{k+1}= &R_{1,3}(\tau \mathcal{A} ) U^k  + \phi_1(\tau \mathcal{A} ) 
 f( U^k,t_k) + \phi_2(\tau \mathcal{A} ) \left( f( a^k,t_k+\tau/2) +f( b^k,t_k+\tau/2) \right)  \\
 &+\phi_3(\tau \mathcal{A} ) f( c^k,t_k+\tau) ,
 \end{split}
 \end{align} 
where
\begin{align*}
a^k&= R_{1,3}(\tau/2 \mathcal{A} )U^k + \phi(\tau/2 \mathcal{A} ) f( U^k,t_k) ,\\
b^k&= R_{1,3}(\tau/2 \mathcal{A} ) U^k +\phi(\tau/2 \mathcal{A} ) f( a^k,t_k+ \tau/2) ,\\
c^k&= R_{1,3}(\tau/2 \mathcal{A} ) a^k +\phi(\tau/2 \mathcal{A} ) 
\left( 2 f( b^k,t_k+ \tau/2) - f( u^k,t_k)\right)  ,
\end{align*} 
with
\begin{align*}
&R_{1,3}(\tau \mathcal{A}) = (24\mathcal{I}- 6 \tau \mathcal{A}) \left( 24\mathcal{I} + 18 \tau \mathcal{A}+6 \tau^2 \mathcal{A}^2 +\tau^3 \mathcal{A}^3\right)^{-1},\\
&R_{1,3}(\tau/2 \mathcal{A}) = 24 (8\mathcal{I} - \tau \mathcal{A}) \left( 192\mathcal{I} + 72 \tau \mathcal{A}+12 \tau^2 \mathcal{A}^2 +\tau^3 \mathcal{A}^3\right)^{-1},\\
&\phi_1(\tau \mathcal{A} )   = \tau (4\mathcal{I} - \tau \mathcal{A}) \left( 24\mathcal{I} + 18 \tau \mathcal{A}+6 \tau^2 \mathcal{A}^2 +\tau^3 \mathcal{A}^3\right)^{-1},\\
&\phi_2(\tau \mathcal{A} )   = 2\tau (4\mathcal{I} + \tau \mathcal{A}) \left( 24\mathcal{I} + 18 \tau \mathcal{A}+6 \tau^2 \mathcal{A}^2 +\tau^3 \mathcal{A}^3\right)^{-1},\\
&\phi_3(\tau \mathcal{A} )   = \tau (4\mathcal{I} + 3\tau \mathcal{A} + \tau^2 \mathcal{A}^2) \left( 24\mathcal{I} + 18 \tau \mathcal{A}+6 \tau^2 \mathcal{A}^2 +\tau^3 \mathcal{A}^3\right)^{-1},\\
&\phi(\tau/2 \mathcal{A} )   = \tau (96\mathcal{I} +12 \tau \mathcal{A} +\tau^2 \mathcal{A}^2 ) \left( 192\mathcal{I} + 72 \tau \mathcal{A}+12 \tau^2 \mathcal{A}^2 +\tau^3 \mathcal{A}^3\right)^{-1},
\end{align*}
where $\mathcal{I}$ is an identity matrix.
\subsection{FFT implementation of ETDRK4-P13 method}

For the efficient implementation of ETDRK4-P13, in the following algorithm the FFT computations of the ETDRK4-P13 method is provided for the case of 3D SFRD with periodic boundary conditions:  

\begin{algorithm}[H]
	\caption{ETDRK4-P13 Method}
	\begin{algorithmic}[1]
		\State Given $ \Lambda^{\alpha/2} $.
		\State  Set $ L=\kappa \tau  \Lambda^{\alpha/2}  $.
		\State  Precompute the following 3D arrays:
		\State $ R_{13} = (24\mathcal{I}  - 6 L ) (24\mathcal{I} + 18 L  +6 L^2 + L^3)^{-1} $.
		\State $ Q_{13} =24 ( 8\mathcal{I} -  L ) (192\mathcal{I} + 72 L + 12 L^2 + L^3)^{-1} $.
		\State $\phi =\tau  (96\mathcal{I}+ 12 L +  L^2) (192\mathcal{I}  + 72 L + 12 L^2 + L^3)^{-1} $.
		\State $\phi_1 =\tau  (4\mathcal{I} - L ) (24\mathcal{I} + 18 L  +6 L^2 + L^3)^{-1} $.
		\State $\phi_2 =2\tau  (4\mathcal{I} + L ) (24\mathcal{I} + 18 L  +6 L^2 + L^3)^{-1} $.
		\State $\phi_3 =\tau  (4\mathcal{I}+3  L +L^2 ) (24\mathcal{I}+ 18 L  +6 L^2 + L^3)^{-1} $.
		\For{$ k=0,1,\cdots, M-1 $} 
		\State $ t_k=k \tau $.
		\State $ u_h=\textit{fftn}(u) $.     \Comment{MATLAB build-in function}
	    \State Compute $ f(u,t_k) $ as $ f_u $.
	    \State  $ a_k= \textit{ifftn}\left( Q_{13} \odot u_h + \phi \odot\textit{fftn}(f_u) \right)  $. \Comment{The operator $ \odot $ stands for element by element multiplication}
	    \State Compute $ f(a_k,t_k+\tau/2) $ as $ f_a $. 
	    \State  $ b_k= \textit{ifftn}\left( Q_{13}\odot  u_h + \phi \odot\textit{fftn}(f_a)  \right)  $.
	    \State Compute $ f(b_k,t_k+\tau/2) $ as $ f_b $. 
	    \State  $ c_k=\textit{ifftn}\left( Q_{13}\odot  \textit{fftn}(a_k)  + \phi \odot \textit{fftn}(2 f_b-f_u) \right) $.
	    \State Compute $ f(c_k,t_k+\tau) $ as $ f_c $. 
	    \State $ u = \textit{ifftn}\left( R_{13} \odot u_h + \phi_1 \odot\textit{fftn}(f_u) + \phi_2 \odot \textit{fftn}(f_a+f_b) + \phi_3 \odot \textit{fftn}(f_c)\right)  $.
		\EndFor

%		\Procedure{Roy}{$a,b$}       \Comment{This is a test}
%		\State System Initialization
%		\State Read the value 
%		\If{$condition = True$}
%		\State Do this
%		\If{$Condition \geq 1$}
%		\State Do that
%		\ElsIf{$Condition \neq 5$}
%		\State Do another
%		\State Do that as well
%		\Else
%		\State Do otherwise
%		\EndIf
%		\EndIf
%		
%		\While{$something \not= 0$}  \Comment{put some comments here}
%		\State $var1 \leftarrow var2$  \Comment{another comment}
%		\State $var3 \leftarrow var4$
%		\EndWhile  \label{roy's loop}
%		\EndProcedure
%		
	\end{algorithmic}
\end{algorithm}
 
 \subsection{Linear analysis}
 The linear truncation error and stability analysis of the ETDRK4-P13  method \eqref{ETDRK4P13} are performed in this subsection.
 \subsubsection{Truncation error of ETDRK4-P13}
 
 It is obvious that the overall spatial discretization is of order four because a fourth-order compact finite difference scheme is applied to fractional Laplacian term.
 Here, we analysis the overall local temporal truncation error for \eqref{ETDRK4P13}  scheme by using Taylor expansion. To this end, we consider the following linear semi-discretization system 
 \begin{align}\label{TE}
 \dfrac{\partial  \mathbf{v} }{\partial  t} + \mathcal{A}  \mathbf{v}  = \mathcal{B}  \mathbf{v} 
 \end{align}
 where $ \mathcal{A} $ and $ \mathcal{B} $ represent matrices derived from the linear spatial discretization of fractional Laplacian and reaction term of a linear space-fractional reaction-diffusion system respectively, and $ \mathbf{v}  $ is a vector of unknowns.

Similar to the previous approach introduced in \cite{zhao2011operator,bhatt2019efficient}, we apply \eqref{ETDRK4P13} to the system \eqref{TE} as follows
 \begin{align}\label{TE1}
 \mathbf{v}^{k+1} = & R_{1,3}(\tau \mathcal{A} ) \mathbf{v}^{k}  + \phi_1(\tau \mathcal{A} ) 
  \mathcal{B} \mathbf{v}^{k} + \phi_2(\tau \mathcal{A} ) \mathcal{B} \left( \mathbf{a}^{k}+\mathbf{b}^{k}\right) + 
 \phi_3(\tau \mathcal{A} ) \mathcal{B} \mathbf{c}^{k},
 \end{align}
 where 
 \begin{align*}
 \mathbf{a}^{k} &= \left(  R_{1,3}(\tau/2 \mathcal{A} ) + \phi(\tau/2 \mathcal{A} ) \mathcal{B} \right) \mathbf{v}^{k} \\
 \mathbf{b}^{k} &= R_{1,3}(\tau/2 \mathcal{A} )   \mathbf{v}^{k}  +  \phi(\tau/2 \mathcal{A} ) \mathcal{B} \mathbf{a}^{k}\\
 \mathbf{c}^{k} &=  R_{1,3}(\tau/2 \mathcal{A} )   \mathbf{v}^{k} +  \phi(\tau/2 \mathcal{A} ) \mathcal{B} \left(2 \mathbf{b}^{k} - \mathbf{v}^{k}  \right) ,
 \end{align*}
 for $ k=0,\cdots,M-1 $.
 By using Taylor expansion, the scheme \eqref{TE1} becomes
 \begin{align}
 \begin{split}
 \mathbf{v}^{k+1} = &\left(  1 + (-\mathcal{A} + \mathcal{B}) \tau + (\mathcal{A}^2 - 2 \mathcal{A} \mathcal{B} + \mathcal{B} ) \tau^2/2 + (-\mathcal{A}^3 +3 \mathcal{A}^2 \mathcal{B} -3 \mathcal{A} \mathcal{B}^2 + \mathcal{B}^3 ) \tau^3/6 +\right.\\
   &\left.   ( \mathcal{A}^4  - 4 \mathcal{A}^3 \mathcal{B} + 6 \mathcal{A}^2 \mathcal{B}^2 - 4 \mathcal{A} \mathcal{B}^3 +\mathcal{B}^4 ) \tau^4/24+ \cdots \right) \mathbf{v}^{k}.
 \end{split}
 \end{align}
For  $ k=0,\cdots,M-1 $, the exact solution of system \eqref{TE} is 
 \[
 \mathbf{v}(t_{k+1}) = e^{(-\mathcal{A} + \mathcal{B} ) \tau}  \mathbf{v}(t_{k}) .
 \]
 Therefore, the local truncation error of ETDRK4-P13  method \eqref{ETDRK4P13} is
 \begin{align}
 \begin{split}
  &\left(  1 + (-\mathcal{A} + \mathcal{B}) \tau + (\mathcal{A}^2 - 2 \mathcal{A} \mathcal{B} + \mathcal{B} ) \tau^2/2 + (-\mathcal{A}^3 +3 \mathcal{A}^2 \mathcal{B} -3 \mathcal{A} \mathcal{B}^2 + \mathcal{B}^3 ) \tau^3/6 +\right.\\
 &\left.   ( \mathcal{A}^4  - 4 \mathcal{A}^3 \mathcal{B} + 6 \mathcal{A}^2 \mathcal{B}^2 - 4 \mathcal{A} \mathcal{B}^3 +\mathcal{B}^4 ) \tau^4/24+ \cdots \right) \mathbf{v}^{k} - 
 e^{(-\mathcal{A} + \mathcal{B} ) \tau} \mathbf{v}(t_{k})  = \mathcal{O}(\tau^5) \mathbf{v}^{k} .
 \end{split}
 \end{align}
 Hence  ETDRK4-P13  method \eqref{ETDRK4P13} is fourth-order in time discretization.
 
 \subsubsection{Stability region of ETDRK4-P13} 
 
 We demonstrate the linear stability of ETDRK4-P13  method \eqref{ETDRK4P13} by plotting its stability region. We consider the following nonlinear autonomous ODE
 \begin{align}
 u_t + c u= f(u). 
 \end{align}
 
 Linearizing the above equation about a fixed point $ g $ such that $ f(g ) - c g =0 $, we obtain the following test equation
 \begin{align}\label{TEq}
  u_t + c u =  \gamma u,
 \end{align}
 where $ u $ is a perturbation of $ g  $ and $ \gamma = f'(g )  $.
 For all $ \gamma $, the fixed point $ g $ is stable if $ \Re (\gamma-c) < 0$.
 
\begin{figure}[H]
  	\centering
  	\includegraphics[width=0.85\linewidth, height=0.3\textheight]{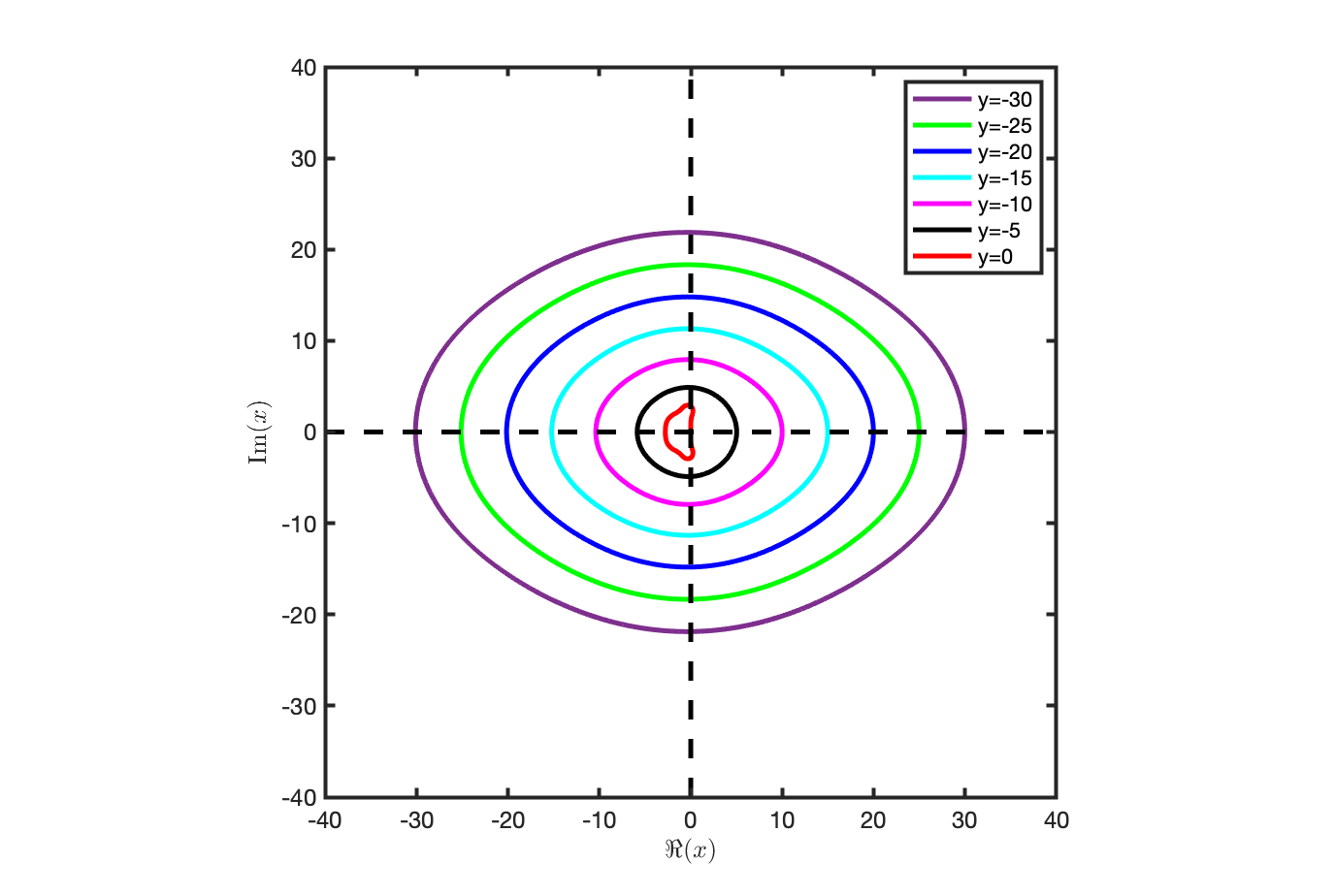}
  	\caption{Stability regions of ETDRK4-P13  method in the complex $ x $-plane.}
  	\label{fig:stp13}
  \end{figure}

 Applying the semi-discrete ETDRK4-P13  method \eqref{ETDRK4P13} to the scalar test equation \eqref{TEq} with $ y=-c \tau $ and $ x=\gamma \tau $, we have the following amplification factor
 
 \begin{align}\label{Eq41}
 r(x,y):=\frac{u^{k+1} }{u^k } = c_0 + c_1 x + c_2 x^2 + c_3 x^3 + c_4 x^4,
 \end{align}
 where 
 \begin{align*}\small
 c_0 = & \frac{6 (4+ y)}{24 - 18 y + 6 y^2 -y^3},\\
 c_1= & \frac{169869312 - 148635648 y - 9621504 y^3 + 470016 y^4 + 165888 y^5 - 
 	46464 y^6 + 6240 y^7 - 600 y^8 + 32 y^9 - y^{10}}{(192 - 72 y + 12 y^2 - y^3)^3(24 - 18 y + 6 y^2 -y^3)},\\
 c_2 =& \frac{442368 - 221184 y + 34560 y^2 + 11712 y^3 - 2976 y^4 + 404 y^5 - 
 	25 y^6 + y^7}{(192 - 72 y + 12 y^2 - y^3)^2(24 - 18 y + 6 y^2 -y^3)},\\
 c_3=& \frac{2 (96 - 12 y + y^2)^2 (1536 - 960 y + 240 y^2 + 8 y^3 + y^4)}{(192 - 72 y + 12 y^2 - y^3)^3(24 - 18 y + 6 y^2 -y^3)},\\
  c_4=&\frac{2 (96 - 12 y + y^2)^3 (4 - 3 y + y^2)}{(192 - 72 y + 12 y^2 - y^3)^3(24 - 18 y + 6 y^2 -y^3)}.
 \end{align*}
 The stability region is four-dimensional and therefore difficult to visualize. To obtain a two- dimensional stability region, we assume $ c $ to be fixed and real and $ \gamma $ is complex.
 The boundaries of the stability regions of the ETDRK4-P13  method are obtained by substituting $ r=e^{i\theta} $, $ \theta\in[0,2\pi] $ into \eqref{Eq41} and solving for $ x $.
 In Fig. \ref{fig:stp13}, we plot the curves of $ | r(x,y)|=1 $ for a complex value of $ x $ and various negative values of $ y $. According to Beylkin et al. \cite{beylkin1998new}, the stability regions should grow as $ y $ approach to $ -\infty $ for the scheme to be useful and practicable. Clearly, the stability regions for the ETDRK4-P13  method as shown in Fig. \ref{fig:stp13} grow larger as $ y\rightarrow  -\infty $ which affirm the stability of the scheme.
  \begin{remark}
  It can be observed from formula \eqref{Eq41} that if $ y=0 $, the amplification factor becomes
  \[
  r(x,0)= 1 + x + \frac12 x^2 + \frac16 x^3 + \frac{1}{24} x^4,
  \]
  and it represents the amplification factor of the fourth-order explicit Runge-Kutta (RK4) scheme.
  \end{remark}

\section{Numerical experiments}
\label{S5}
This section discusses the results obtained by applying the proposed method to various widely known SFRDEs including 2D FitzHugh-Nagumo, Gierer-Meinhardt, Gray-Scott and 3D Schnakenberg models.
% Numerical experiments include Fisher's and Huxley equations with a source term, FitzHugh Nagumo model which arises in the excitable media study [cite], Gray–Scott model and finally Schnakenberg model that describes an auto-catalytic chemical reaction [cite].
The accuracy of the scheme is measured in terms of maximum error norm $\norm{\cdot}_\infty $ defined as 
\[
\| u-u_{h,\tau} \|_{\infty} = \max_{n} | u(\mathbf{x}_n,T) - u_n^M|
\]
where $ u(\mathbf{x}_n,T) $ and $ u_n^M $ are the $ n $th exact and numerical solution at the final time $ T $, respectively. In additional, the order of convergence in space and time is computed as 
\[
Order =\log_2(E_{h,\tau}/E_{h/2,\tau/2}) ,
\]
 where $  E_{h,\tau} = \| u-u_{h,\tau} \|_{\infty} $ and $  E_{h/2,\tau/2} = \| u-u_{h/2,\tau/2} \|_{\infty} $ with spatial step size $ h $ and temporal size $ \tau $.
The numerical experiments are performed in MATLAB R2019b platform on a MacBook Pro with 3.1 GHz Dual-Core Intel Core i7 CPU and 8 GB memory.
 
\subsection{Stability, accuracy and efficiency tests}
The following examples are considered as a benchmark problems in order to investigate the performance in terms of accuracy and efficiency of the proposed method.

{\bf Example 1.}~~	In this example, we consider the following 1D space-fractional Fisher's equation with source term over a domain $ \Omega=[0,1] $
	\begin{align}\label{Ex1}
	\dfrac{\partial  u }{\partial  t} &= - \kappa (-\Delta)^{\alpha/2} u +  f(u,t), 
	\end{align}
	where
	\[
	 f(u,t) = - 2e^{-t} \sin^3(2\pi x)+\frac{\kappa}{4}\kappa e^{-t} \left( 3 (2 \pi)^{\alpha} \sin(2\pi x) - (6 \pi)^{\alpha} \sin(6\pi x) \right) + e^{-2t} \sin^6(2\pi x) + u - u^2,
	\]
	subject to homogenous Dirichlet boundary conditions and the following initial condition
	\[
	g(x)= \sin^3(2\pi x), \quad 0\leq x \leq 1.
	\]
The exact solution to \eqref{Ex1} is 
\begin{equation*}
u(x,t)= e^{-t} \sin^3(2\pi x).
\end{equation*}

\begin{figure}[H]
	\centering
	\includegraphics[width=1.1\linewidth, height=0.3\textheight]{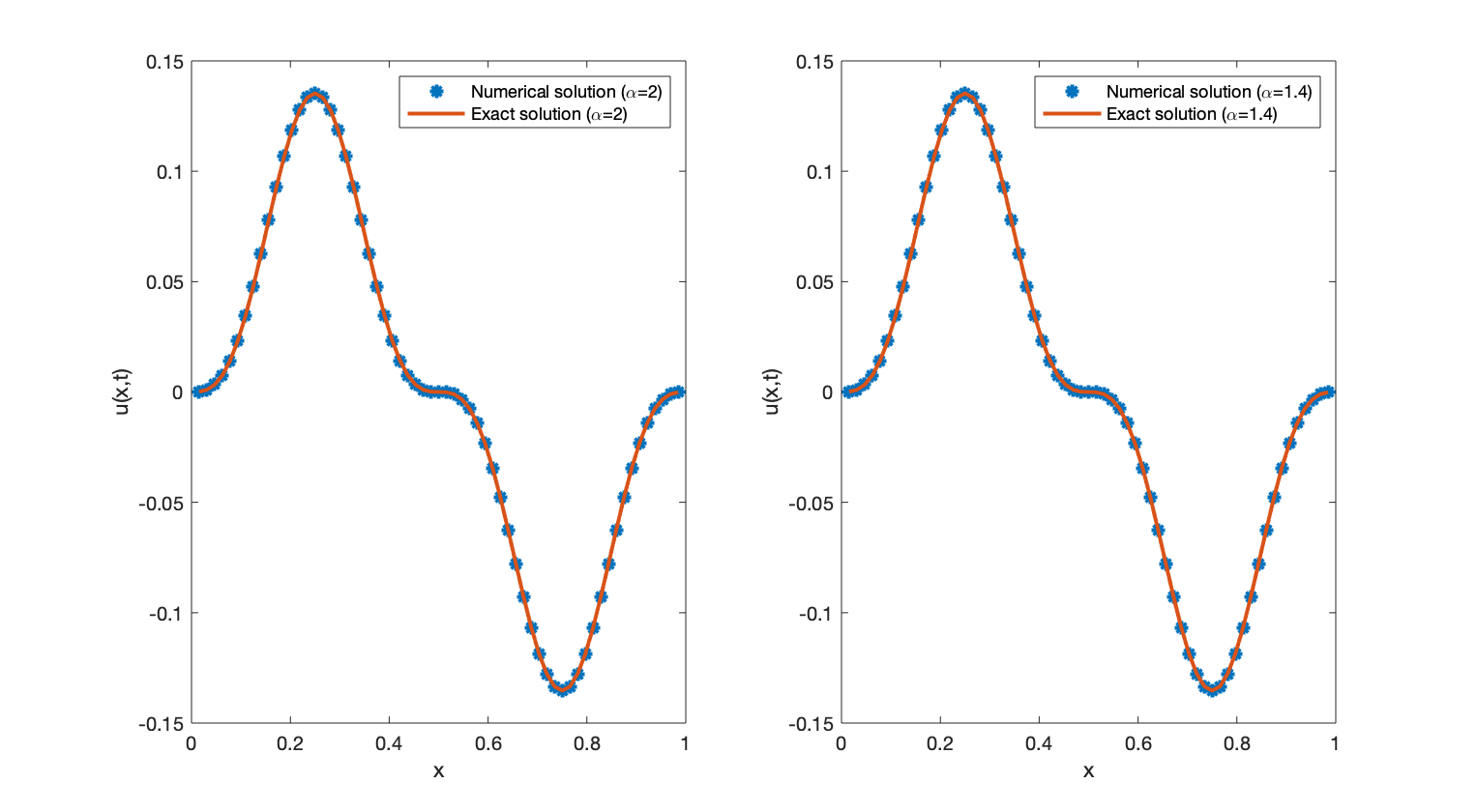}
	\caption{Numerical solutions of Example 1 vs. the exact solution. }
	\label{fig:ex1}
\end{figure}

In the first set of experiments, we compared the analytical solution with solution obtained by the proposed method in order to see how the proposed method is able to capture the spatio temporal solution profile of the component $u$. Fig. \ref{fig:ex1} displays the approximate solution obtained via the proposed method up to $ T=2 $ and the corresponding exact solution for $ \alpha=2 $ and $ \alpha=1.4 $ with $ \kappa=10 $, $ N=64 $ and $ \tau=0.004 $. It is clear from the Fig. \ref{fig:ex1} that the numerical solution provided by the proposed method is in good agreement with an exact solution.

\begin{figure}[H]
	\centering
	\includegraphics[width=1.1\linewidth, height=0.3\textheight]{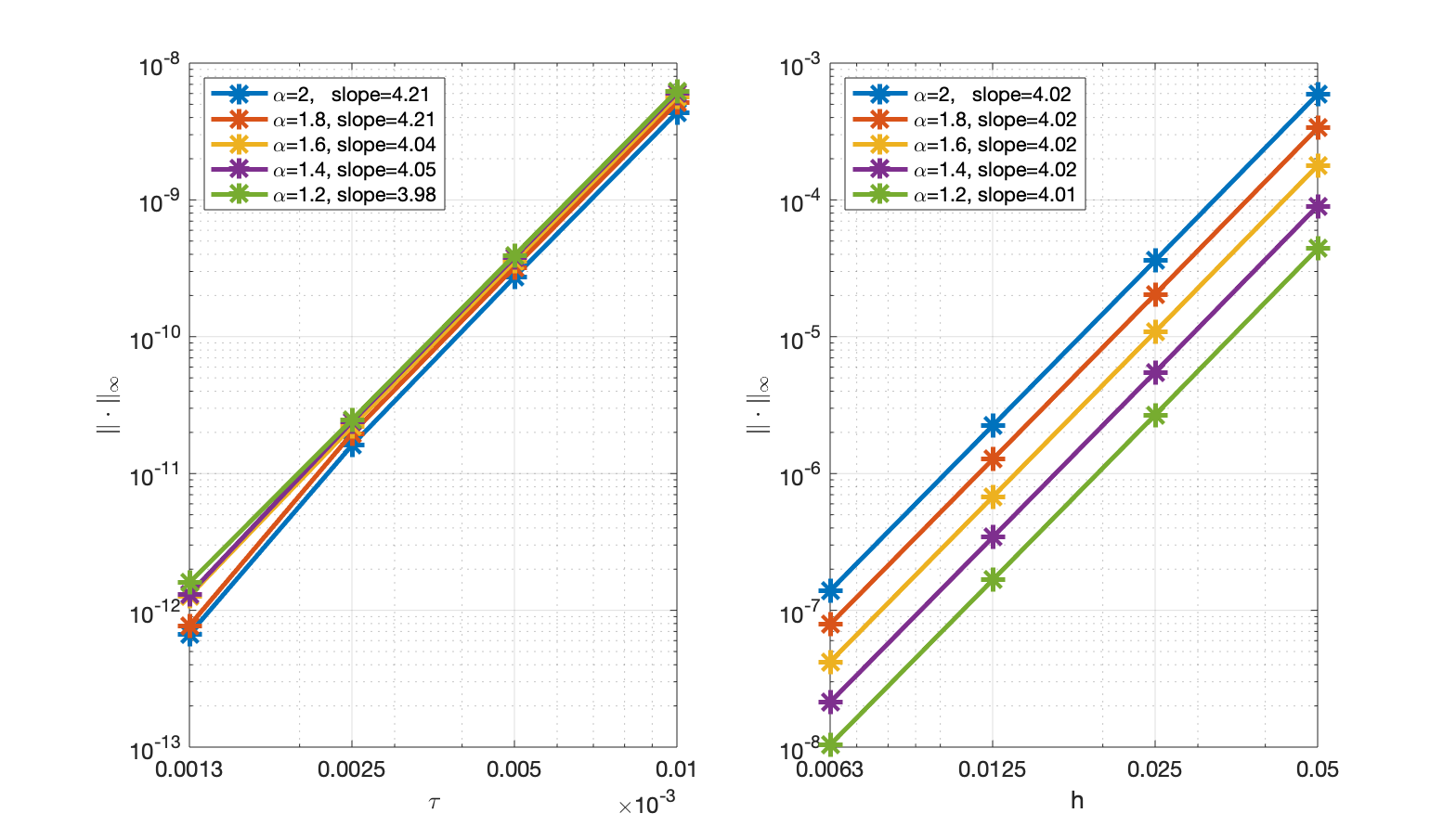}
	\caption{Left label: The log-log graph of errors vs. temporal step size $ \tau $ with $ N=3000 $ and $ \kappa=1e-3 $. Right label: The log-log graph  of errors vs. spatial step size $ h $ with $ M=1000 $ and $ \kappa=1e-3 $. }
	\label{fig:ex1slop}
\end{figure}

In order to do the empirical convergence of the proposed method, we ran a another set of experiments on Example 1 at $T=1$ and captured the maximum errors for various values of $\alpha,~h$ and $k$ with $\kappa=1$. In order to visualize the temporal and space rates of convergence of the proposed method, we depicted them in Fig. \ref{fig:ex1slop} with log-log scale graph. One can notice that the computed rates of convergences agree with expected order of convergence of the scheme which is a fourth-order accurate in both space and time.

A third set of experiments ran on Example 1 in order to compare the performance of the proposed method in terms of accuracy and efficiency with the results obtained via the RK4 method for $ \alpha=1.8 $ and different values of $ h $ and $ \tau $. Maximum errors, rates of convergence, and CPU times for the  methods are reported in Table \ref{table:comp}. From Table \ref{table:comp}, one can see that the temporal step size must be proportional to $\alpha$ power of the spatial step size $h$ as expected when the explicit RK4 is used due to the stability constraint from the fractional diffusion term. From Table \ref{table:comp}, one can see that the computational cost for the RK4 method is much higher than the proposed method at the same spatial step size $h$ with similar size of errors and order of accuracy. So, the results in Table \ref{table:comp} justify well enough abandoning the implementation of the RK4 method in the remaining part of this work.
\begin{table}[H]                                                                     
	\centering     
	\caption{Comparison of accuracy and efficiency of the ETDRK4-P13 and RK4 methods for Example 1 \ref{Ex1} with $ \kappa=10 $.}                                                          
	\label{table:comp}                                                                        
	\begin{tabular}{cccccccccc}    
		\toprule
		\multicolumn{1}{c}{} &\multicolumn{2}{l}{ETDRK4-P13} & &\multicolumn{3}{l}{RK4} &\multicolumn{3}{l}{RK4}  \\
		\multicolumn{1}{c}{} &\multicolumn{2}{l}{$\alpha=1.8$, $\tau =\frac{h}{4 \kappa}  $} & &\multicolumn{3}{l}{$\alpha=1.8$, $\tau =\frac{h}{4 \kappa}  $} &\multicolumn{3}{l}{$\alpha=1.8$, $\tau =\frac{h^{\alpha}}{4 \kappa}  $}  \\
		\cmidrule(r){2-4}  \cmidrule(r){5-7} \cmidrule(r){8-10}
		$ h $&$ \| \cdot \|_{\infty}  $&Order&CPU(s)& $ \| \cdot \|_{\infty} $&Order&CPU(s)&$ \| \cdot \|_{\infty} $&Order&CPU(s)\\                                        
		\midrule                                                                                         
		$ \frac{1}{8} $  & 1.3871e-02 & - & 0.0404 & NaN & - & 0.0332 & 1.3865e-02 & - & 0.2278 \\
		                                                                        
		$ \frac{1}{16} $ & 7.2947e-04 & 4.25 & 0.0617 & NaN & - & 0.0757 & 7.5050e-04 & 4.21 & 0.6722 \\  
		                                                                           
		$ \frac{1}{32} $& 4.3772e-05 & 4.06 & 0.1440 & NaN & - & 0.1712 & 4.3772e-05 & 4.10 & 2.7245 \\  
		                                                                      
		$ \frac{1}{64} $  & 2.7084e-06 & 4.01 & 0.3625 & NaN & - & 0.4089 & 2.9474e-06 & 3.89 & 11.7495 \\ 
			\bottomrule                                                                              
\end{tabular}                                                                                                                        
\end{table}       
%%%%%%%%%%%%%

%\subsection{Convergence and efficiency tests}
{\bf Example 2.}~~Consider the following 2D space-fractional Huxley equation with source term over a domain $ \Omega=[0,1]^2 $
	\begin{align}\label{Ex2}
	\dfrac{\partial  u(x,y,t) }{\partial  t} &= - \kappa (-\Delta)^{\alpha/2} u + u (1-u) (u-1)+ f(x,y,t), 
	\end{align}
	where
	\begin{align*}
	f(x,y,t) =& \alpha t^{\alpha-1} \cos^3(2\pi x) \cos^3(2\pi y) + t^{\alpha} \varPhi(x,y) -  t^{\alpha} \cos^3(2\pi x) \cos^3(2\pi y) \left(1-t^{\alpha} \cos^3(2\pi x) \cos^3(2\pi y) \right) \\
	& \times \left( t^{\alpha} \cos^3(2\pi x) \cos^3(2\pi y) -1\right) ,
	\end{align*}
with
\begin{align*}
\varPhi(x,y) = &\frac{\kappa}{16} \left( 9 (8\pi^2 )^{\alpha/2} \cos(2\pi x) \cos(2\pi y) + 3 (40\pi^2 )^{\alpha/2} \cos(6\pi x) \cos(2\pi y) +3 (40\pi )^{\alpha/2} \cos(2\pi x) \cos(6\pi y) \right.\\
 &\left. + (72\pi^2 )^{\alpha/2} \cos(6\pi x) \cos(6\pi y) \right) ,
\end{align*}
	subject to periodic and homogenous Neumann boundary conditions. The initial condition can be obtained from the analytical solution $ u(x,y,t)=t^{\alpha} \cos^3(2\pi x)\cos^3(2\pi y) $.

To investigate the order of accuracy of the method for varying fractional order $\alpha$, a numerical test on Example 2 was performed. In the computation, a simulation was run up to the final time $ T=1 $ with $ \kappa=1 $ by initially setting $h=0.1$ and $\tau=0.1h$ then reduced both of them by a factor of 2 in each refinement. The maximum error and rates of convergence for $\alpha = 1.2,1.4,1.8,2.0$ of the diffusion process are listed in Table \ref{table:Ex2}. As expected we observe from Table \ref{table:Ex2} that the computed convergence rates of the proposed method apparently demonstrate the expected fourth-order accuracy in both time and space for varying $\alpha$. 

%Numerical results at the final time $ T=1 $ with $ \kappa=1 $ for Example 2 are reported in Table \ref{table:Ex2}. As expected, we again observe a fourth-order accuracy of both space and time discretization for proposed method. Moreover, one can observe that all CPU times are quite small due to the implementation of FFT.

\begin{table}[H]                                                           
	\centering   
	\caption{Numerical errors, convergence rates and the CPU times for Example 2 with homogenous Neumann or periodic boundary conditions $ (\tau=0.1 h)$.}                                                
	\label{table:Ex2}                                                                   
	\begin{tabular}{ccccccccc}                                      
		\toprule
		
		&&\multicolumn{3}{c}{Neumann boundary condition} & &\multicolumn{3}{c}{Periodic boundary condition}  \\
			\cmidrule(r){3-5}  \cmidrule(r){7-9} 
		$ \alpha $	&$ 1/h  $&  $ \| \cdot \|_{\infty} $ & Order & CPU(s) &&$ \| \cdot \|_{\infty} $ & Order & CPU(s)\\
		
		 \midrule                                                             
		2.0 & 10 & 2.5214e-02 & -& 0.0564 && 2.5214e-02 & - & 0.0317 \\   
		                                                               
		 & 20 & 1.4159e-03 & 4.15 & 0.2209 && 1.4159e-03 & 4.15 & 0.1897 \\   
	                                                              
		 & 40 & 8.6173e-05 & 4.04 & 1.2567 && 8.6173e-05 & 4.04 & 1.0425 \\   
	                                                            
		 & 80 & 5.3493e-06 & 4.01 & 8.5284 && 5.3493e-06 & 4.01 & 7.6987 \\   
	                                                           
		 & 160 & 3.3373e-07 & 4.00 & 66.3996 && 3.3373e-07 & 4.00 & 58.2910 \\
		&&&&&&&&\\                                                          
		1.8 & 10 & 2.2581e-02 & - & 0.0583 && 2.2581e-02 & - & 0.0307 \\   
		                                                            
		 & 20 & 1.2724e-03 & 4.15 & 0.1989 && 1.2724e-03 & 4.15 & 0.1779 \\   
		                                                                 
		 & 40 & 7.7475e-05 & 4.04 & 1.2673 && 7.7475e-05 & 4.04 & 0.9911 \\   
		                                                              
		 & 80 & 4.8107e-06 & 4.01 & 8.4887 && 4.8107e-06 & 4.01 & 7.7161 \\   
		                                                              
		 & 160 & 3.0017e-07 & 4.00 & 66.3678 && 3.0017e-07 & 4.00 & 58.1880 \\
		&&&&&&&&\\                                                              
		1.4 & 10 & 1.7272e-02 & - & 0.0498 && 1.7272e-02 & - & 0.0376 \\   
		                                                          
		 & 20 & 9.8009e-04 & 4.14 & 0.2529 && 9.8009e-04 & 4.14 & 0.1348 \\   
		                                                             
		 & 40 & 5.9702e-05 & 4.04 & 1.2352 && 5.9702e-05 & 4.04 & 0.9927 \\   
		                                                              
		 & 80 & 3.7069e-06 & 4.01 & 8.5357 && 3.7069e-06 & 4.01 & 7.7176 \\   
		                                                                 
		 & 160 & 2.3115e-07 & 4.00 & 68.6389 && 2.3115e-07 & 4.00 & 59.5034 \\
		&&&&&&&&\\                                                                 
		1.2 & 10 & 1.4604e-02 & - & 0.0497 && 1.4604e-02 & - & 0.0304 \\   
		                                                                
		 & 20 & 8.3194e-04 & 4.13 & 0.2907 && 8.3194e-04 & 4.13 & 0.1361 \\   
		                                                                
		 & 40 & 5.0661e-05 & 4.04 & 1.2462 && 5.0661e-05 & 4.04 & 0.9914 \\   
		                                                              
		 & 80 & 3.1332e-06 & 4.02 & 8.5007 && 3.1332e-06 & 4.02 & 7.7329 \\   
		                                                              
		 & 160 & 1.9001e-07 & 4.04 & 67.0404 && 1.9001e-07 & 4.04 & 58.6672 \\
		\bottomrule                                                                    
	\end{tabular}                                                           
	                                      
\end{table}     

\subsection{Some applications}

\subsubsection{\textbf{2D Fitzhugh-Nagumo model}: }

The study of the excitable media is made possible using FitzHugh-Nagumo (FHN) model. A cubic nonlinear reaction term is used to model the propagation of the transmembrane potential in the nerve axon  \cite{fitzhugh1961impulses,nagumo1962active}
\begin{align}
\begin{split}
\dfrac{\partial  u }{\partial  t} &= - \kappa (-\Delta)^{\alpha/2} u + u (1-u) (u-\mu) - v,\\
\dfrac{\partial  v }{\partial  t} &= \epsilon \left( \beta u - \gamma v -\delta \right).
\end{split}
\end{align}
The interest domain is $ \Omega=[0,2.5]^2 $ and the initial conditions are taken as
\begin{align*}
\begin{split}
&u(x,y,0)= 
\begin{cases}
&1, \hspace{1cm} 0 < x \leq 0.125,\, 0 < y < 0.125, \\
& 0, \hspace{1cm}  \text{elsewhere},
\end{cases}\\
&v(x,y,0)= 
\begin{cases}
&0.1, \hspace{1cm} 0 < x < 2.5,\, 0.125 \leq  y < 2.5, \\
& 0, \hspace{1.2cm}  0 < x < 2.5, \,0 < y <0.125.
\end{cases}
\end{split}
\end{align*}
The choice of parameters $ \mu=0.1 $, $ \epsilon=0.01 $, $ \beta=0.5 $, $ \gamma=1 $ and $ \delta=0 $ is known to generate stable patterns in the system in the form of spiral waves.
% where $ \mu=0.1 $, $ \epsilon=0.01 $, $ \beta=0.5 $, $ \gamma=1 $ and $ \delta=0 $.

In order to illustrate the effect of varying $ \alpha $ and $ \kappa $, the model is considered subject to periodic \cite{zhang2020stabilized} or homogenous Neumann boundary conditions \cite{bueno2014fourier} with $ N=256 $ and $ \tau=1 $. 

\begin{figure}[H]
	\centering
	\includegraphics[height=3in,width=5in]{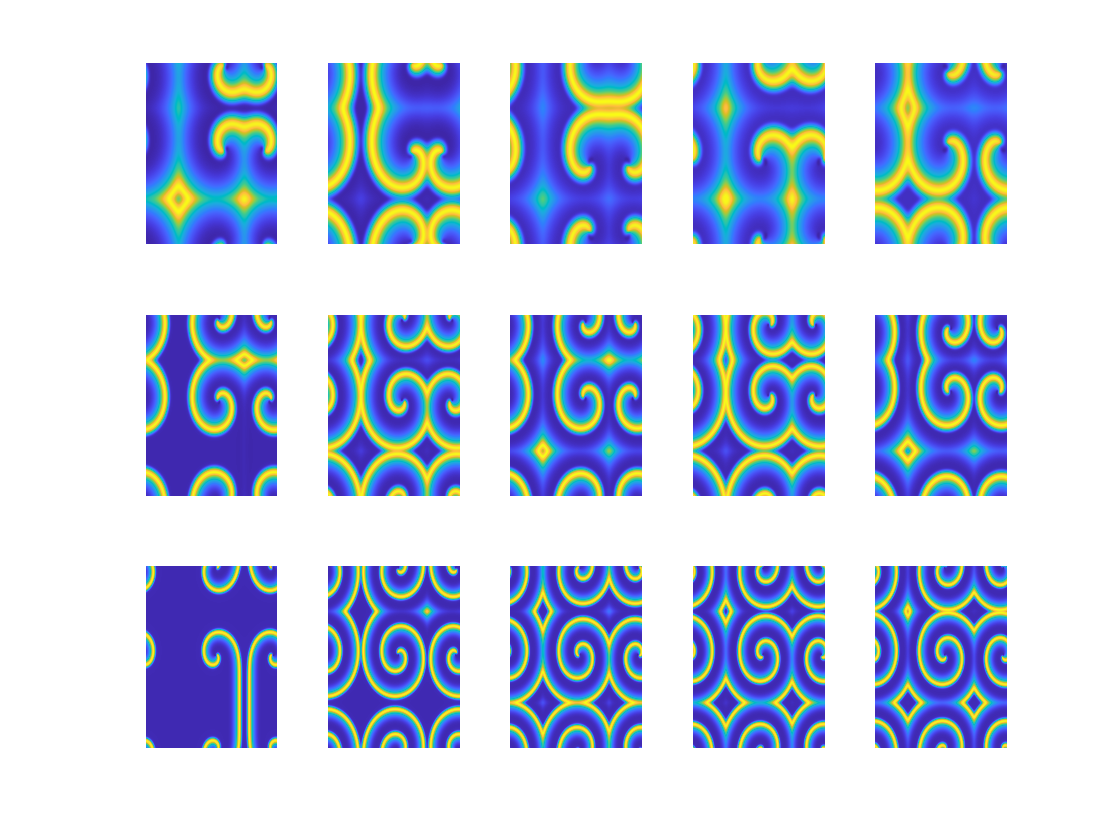}
	\caption{Evolution of the solution $ v $ with $ \kappa=1e-4 $ for varying $ \alpha=2,1.7,1.5$ (Up-Down) at different time $ T=400,700,1000,1500,2000 $ (Left-Right).}
	\label{fig:fhnp}
\end{figure}

\begin{figure}[H]
	\centering
	\includegraphics[height=1in,width=4in]{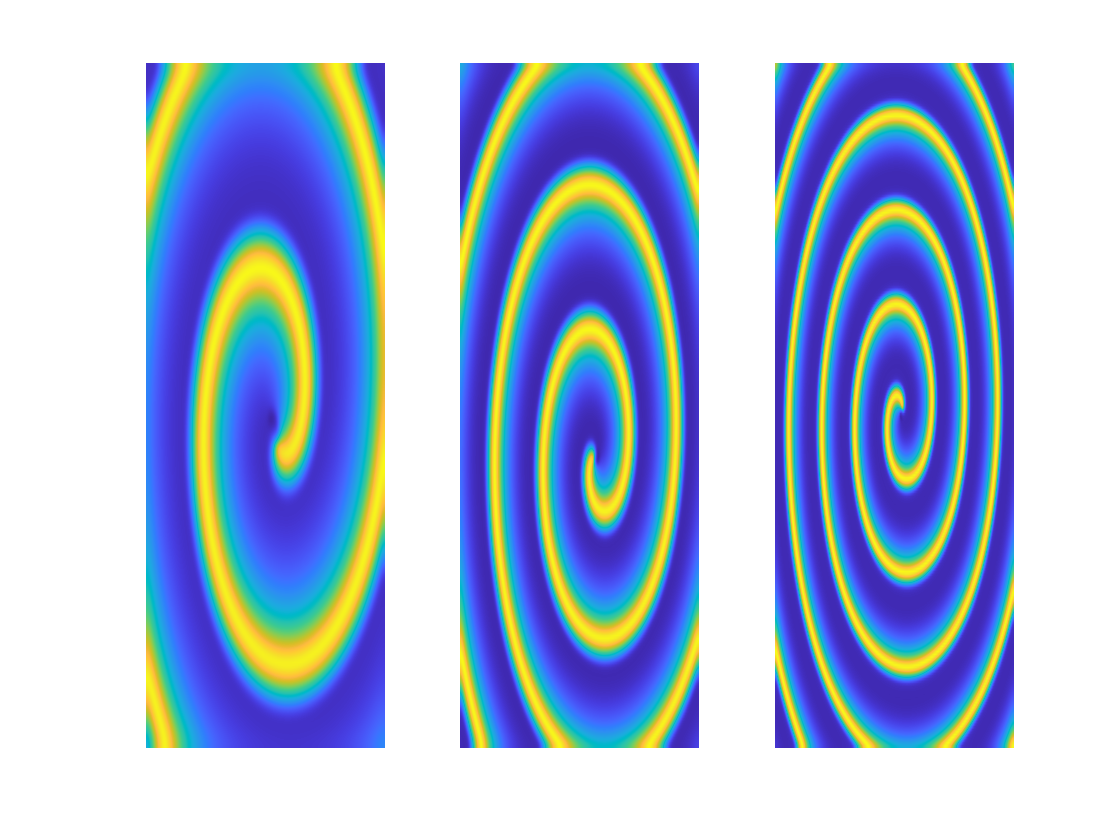}
	\caption{Spiral waves for varying $ \alpha=2, 1.7, 1.5 $ (Left-Right) and $ \kappa=1e-4 $ at $ T=2000 $.}
	\label{fig:fhnk}
\end{figure}

\begin{figure}[H]
	\centering
	\includegraphics[height=1in,width=4in]{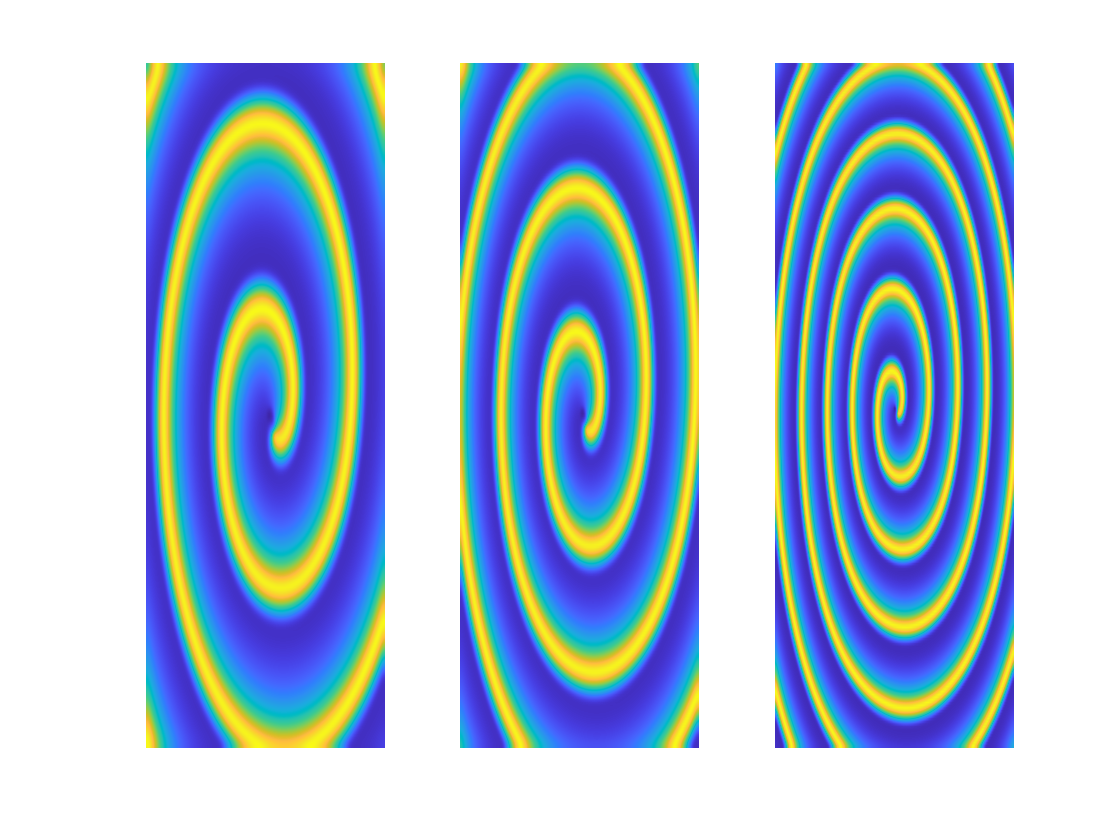}
	\caption{Spiral waves for $ \alpha=2 $ and varying diffusion coefficient $ \kappa=5e-5, 3e-5, 1e-5 $ (Left-Right) at $ T=2000 $.}
	\label{fig:fhna2}
\end{figure}

The aerial views of the concentration profile of the component $v$ subject to periodic boundary conditions with $ \kappa=1e-4 $ at different times and various $\alpha$ values are presented in Fig. \ref{fig:fhnp} to illustrate the effect of fractional diffusion  in FHN model.  From Fig. \ref{fig:fhnp}, it can be clearly depicted that the width of the wavefront is reduced for decreasing fractional power $\alpha$. 
%Figure. \ref{fig:fhnp} presents the stable rotating solution at different time to demonstrate the effect of fractional order in the model subject to periodic boundary conditions with $ \kappa=1e-4 $. It is remarkable that the width of excitation wavefront is reduced when fractional power $ \alpha $ is decreasing.
In the second set of numerical experiments, a stable rotating solution of the component $v$ at $ T=2000 $ is also presented in Fig. \ref{fig:fhnk} to illustrate the effect fractional diffusion and diffusion coefficient in the model subject to homogenous Neumann boundary conditions with $ \kappa=1e-4 $ and varying $ \alpha $. 
It can be clearly seen from Fig. \ref{fig:fhnk} that the role of reducing the fractional power $ \alpha $ is not equivalent to the influence of a decreased diffusion coefficient in the pure diffusion case (Fig. \ref{fig:fhna2}).

\subsubsection{\textbf{2D Gierer-Meinhardt model}: }
\begin{figure}[H]
	\centering
	\includegraphics[height=3in,width=6in]{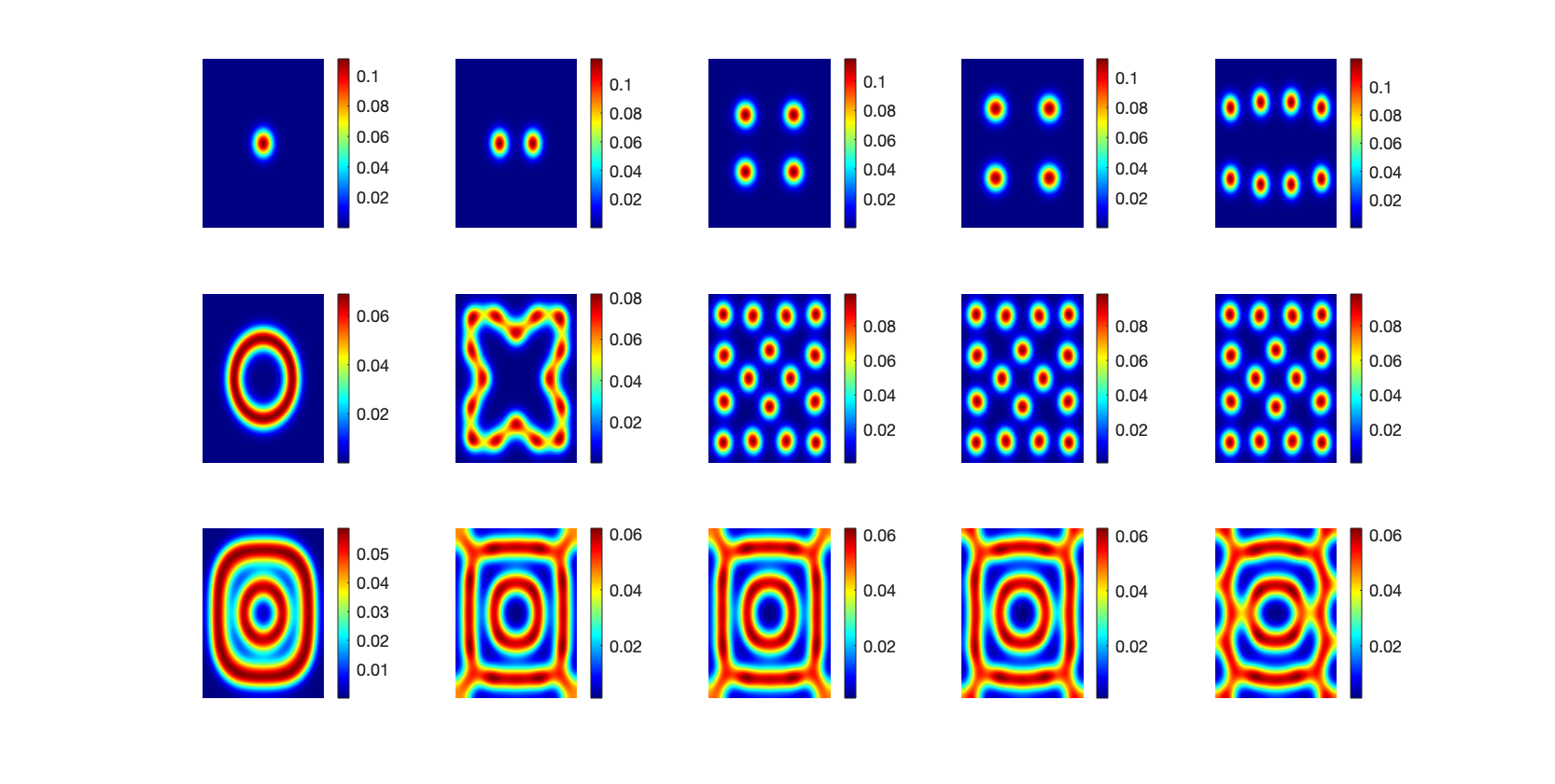}
	\caption{Turing pattern in the Gierer-Meinhardt model \eqref{GM} for $ K=0.0162 $, $ \alpha=2 $ and varying $ \beta=2, 1.9, 1.8$ (Up-Down) at $ T=100, 300, 500, 700, 1000 $ (Left-Right).}
	\label{fig:gmk16}
\end{figure}

\begin{figure}[H]
	\centering
	\includegraphics[height=3in,width=6in]{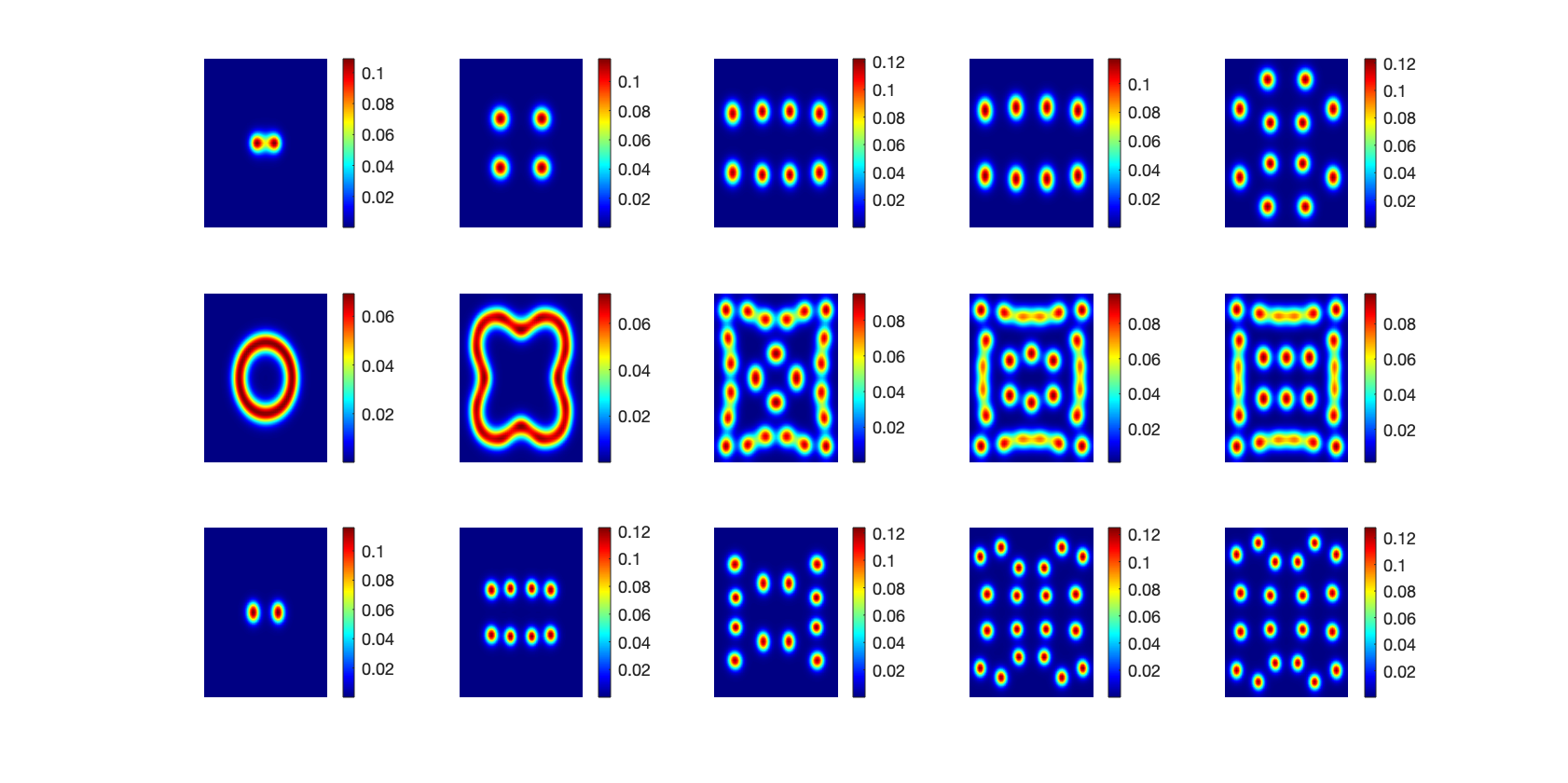}
	\caption{Turing pattern in the Gierer-Meinhardt model \eqref{GM} for $ K=0.0162 $ and varying $ (\alpha, \beta)=(1.9,1.9), (1.9,1.8), (1.8,1.8) $ (Up-Down) at $ T=100, 300, 500, 700, 1000 $ (Left-Right).}
	\label{fig:gmk16a}
\end{figure}

\begin{figure}[H]
	\centering
	\includegraphics[height=3in,width=6in]{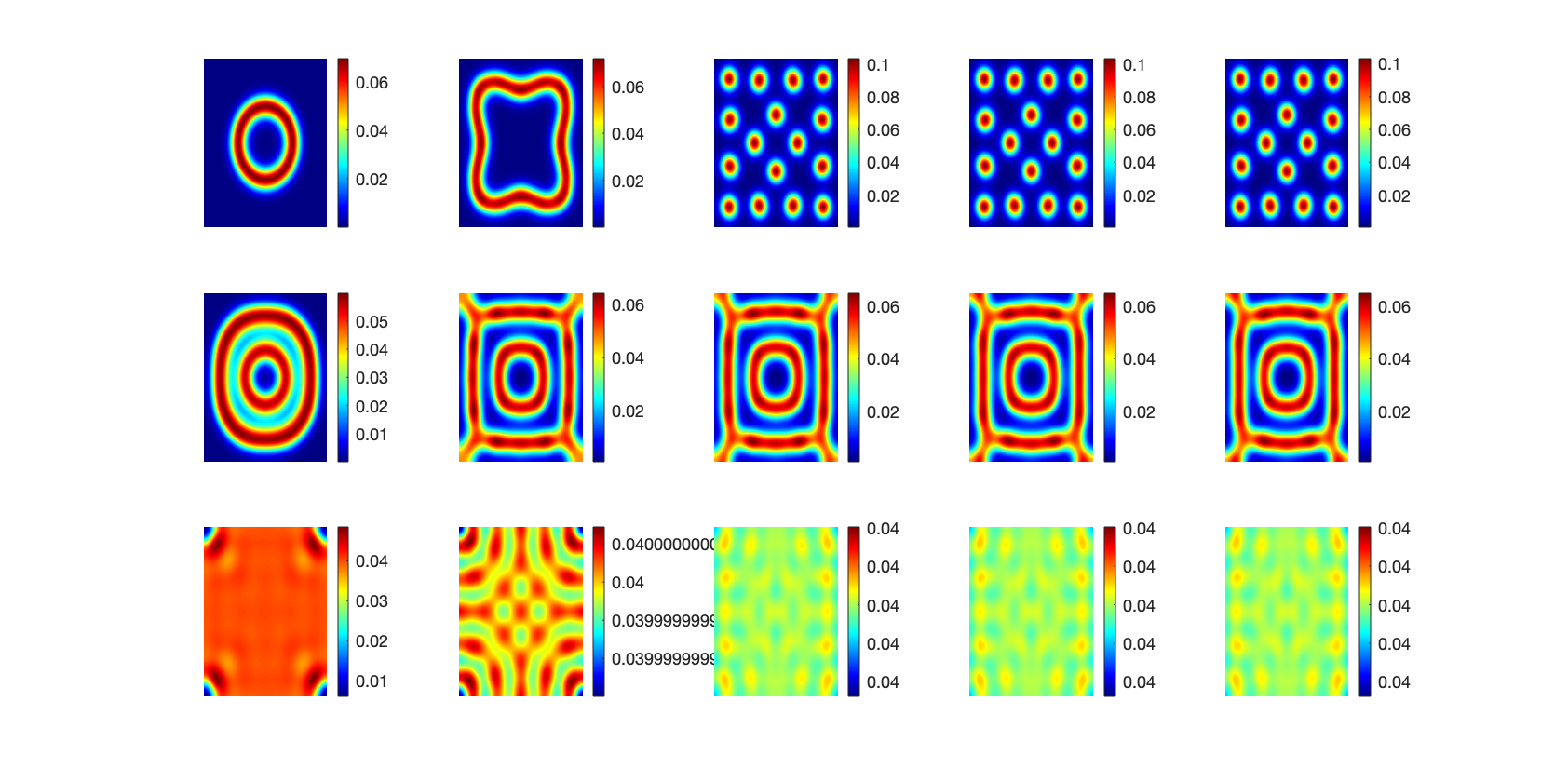}
	\caption{Turing pattern in the Gierer-Meinhardt model \eqref{GM} for $ K=0.0128 $, $ \alpha=2 $ and varying $ \beta=2, 1.9, 1.8$ (Up-Down) at $ T=100, 300, 500, 700, 1000 $ (Left-Right).}
	\label{fig:gmk13}
\end{figure}

\begin{figure}[H]
	\centering
	\includegraphics[height=3in,width=6in]{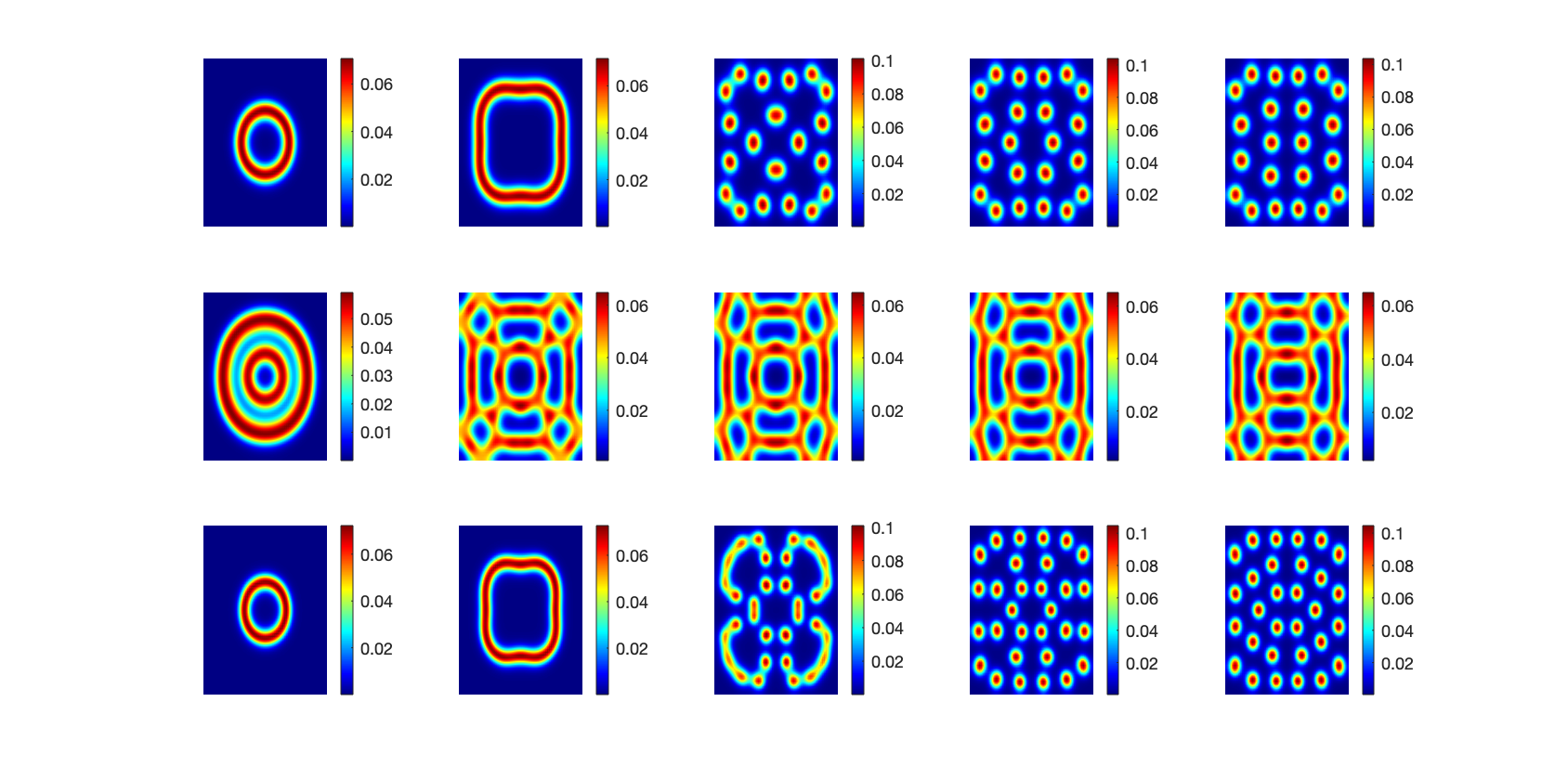}
	\caption{Turing pattern in the Gierer-Meinhardt model \eqref{GM} for $ K=0.0128 $ and varying $ (\alpha, \beta)=(1.9,1.9), (1.9,1.8), (1.8,1.8) $ (Up-Down) at $ T=100, 300, 500, 700, 1000 $ (Left-Right).}
	\label{fig:gmk13a}
\end{figure}

The Gierer–Meinhardt model is widely used in the study of some basic phenomena in the process of morphogenesis. Here, we consider the following fractional Gierer–Meinhardt model \cite{gierer1972theory}:
\begin{align}\label{GM}
\begin{split}
\dfrac{\partial  u }{\partial  t} &= - \kappa_u (-\Delta)^{\alpha/2} u + \frac{u^2}{v} -u \\
\dfrac{\partial  v }{\partial  t} &= - \kappa_v (-\Delta)^{\beta/2} v + \frac{u^2}{\epsilon \mu}  -\frac{v}{\mu} ,
\end{split}
\end{align}
on $ \Omega=[-1,1]^2 $ with homogenous Neumann boundary condition. The parameters are chosen as $ \epsilon=0.04 $, $ \mu=0.1 $. The diffusion coefficient are taken as $ \kappa_u=\epsilon^2 $ and $ \kappa_v= \frac{K}{\mu}$.
The initial condition are chosen as 
\begin{align*}
& u(x,y,0) = \frac12 \left( 1+ 0.001 \sum_{j=1}^{20} \cos(\frac{\pi}{2} j y) \right) {\sech}^2\left( \frac{\sqrt{x^2+y^2}}{2 \epsilon}\right) ,\\
& v(x,y,0) =  \frac{\cosh(1-\sqrt{x^2+y^2})}{3 \cosh(1)}.
\end{align*}
In this set of numerical experiments, we chose $ N=64, \tau=0.1$ and $ K=0.0162 $ and ran the simulation profile of the component $ u $ for fixed $ \alpha=2 $ and different values of $ \beta $ at different times $ T $ and depicted in Fig. \ref{fig:gmk16}. As shown in \cite{wang2019turing}, Turning instability will appear for $ \beta < 1.8 $, so the results were omitted here. In Fig. \ref{fig:gmk16}, it can be observed that the stripe pattern for $ \beta=1.8 $. If we keep increasing the value of $ \beta $, then a chain cluster of Tuning spotted patterns initially evolve. The pure Turing spot patterns are clearly obtained when the simulation time increases to $ 1000 $; as seen in Fig. \ref{fig:gmk16}. In Fig. \ref{fig:gmk16a}, we plot the solution of $ u $ for different $\alpha$ and $\beta$ at different times $ T $. Again, a chain cluster of Turning spotted patterns appear as time increases.
In the second set of numerical experiments, we consider the case for $ K=0.0128 $. Figs. \ref{fig:gmk13} and \ref{fig:gmk13a} exhibit the scenarios of pattern formation for different values of $\alpha$ and $\beta$ at different times $ T $. In Fig. \ref{fig:gmk13}, we fix $ \alpha=2 $ and vary $ \beta $. When $ \beta=1.8 $, then the system become stable, that is, the solution approaches the steady state $ \epsilon=0.04 $ while Turing patterns appear when $ \beta $ increases, as seen in Fig. \ref{fig:gmk13}. In Fig. \ref{fig:gmk13a}, one can observe that a stripe pattern turns into a spot pattern for different $\alpha$ and $\beta$. For more details, we refer to  \cite{wang2019turing,fernandes2012adi,qiao2008numerical,mccourt2008spectral}.

\subsubsection{\textbf{2D Gray-Scott model}: }
The fractional Gray-Scott model that describes an autocatalytic reaction-diffusion process between two chemical species with concentrations $ u $ and $ v $ is written as \cite{gray1983autocatalytic,gray1985sustained}

\begin{align}\label{SG}
\begin{split}
\dfrac{\partial  u }{\partial  t} &= - \kappa_u (-\Delta)^{\alpha/2} u -u v^2 + F (1-u),\\
\dfrac{\partial  v }{\partial  t} &= - \kappa_v (-\Delta)^{\alpha/2} v +u v^2 - (F+K) v,
\end{split}
\end{align}
The diffusion rates in the process satisfy $ \kappa_u\geq 0 $, $ \kappa_v\geq 0 $. The parameters $ F $, $ K $ are the dimensionless feed rate and decay rate, respectively.
The spatially initial condition is
\begin{align*}
(u,v)=
\begin{cases}
& (1,0), \hspace{.3cm} (x,y) \in \Omega\setminus O ,\\
&(\frac12,\frac14), \hspace{.3cm} (x,y) \in O,
\end{cases}
\end{align*}
where $ \Omega=(0,1)^2 $ and $ O = \{ (x,y) | (x-0.5)^2 + (y-0.5)^2 \leq 0.0016 \} $.
We take $ N=265 $, $ \tau=1 $, $ F=0.03 $ and varying $ K $, $ \alpha $. The ratio of diffusion coefficients are chosen as $ \kappa_u=2e-5 $, $ \kappa_v=1e-5 $ $ (\kappa_u/\kappa_v > 1) $ to make the model generating different types of pattern formation \cite{wang2019fractional}.

\begin{figure}[H]
	%\centering
	\begin{subfigure}[h]{1\textwidth}
		\includegraphics[height=2.2in,width=5.2in]{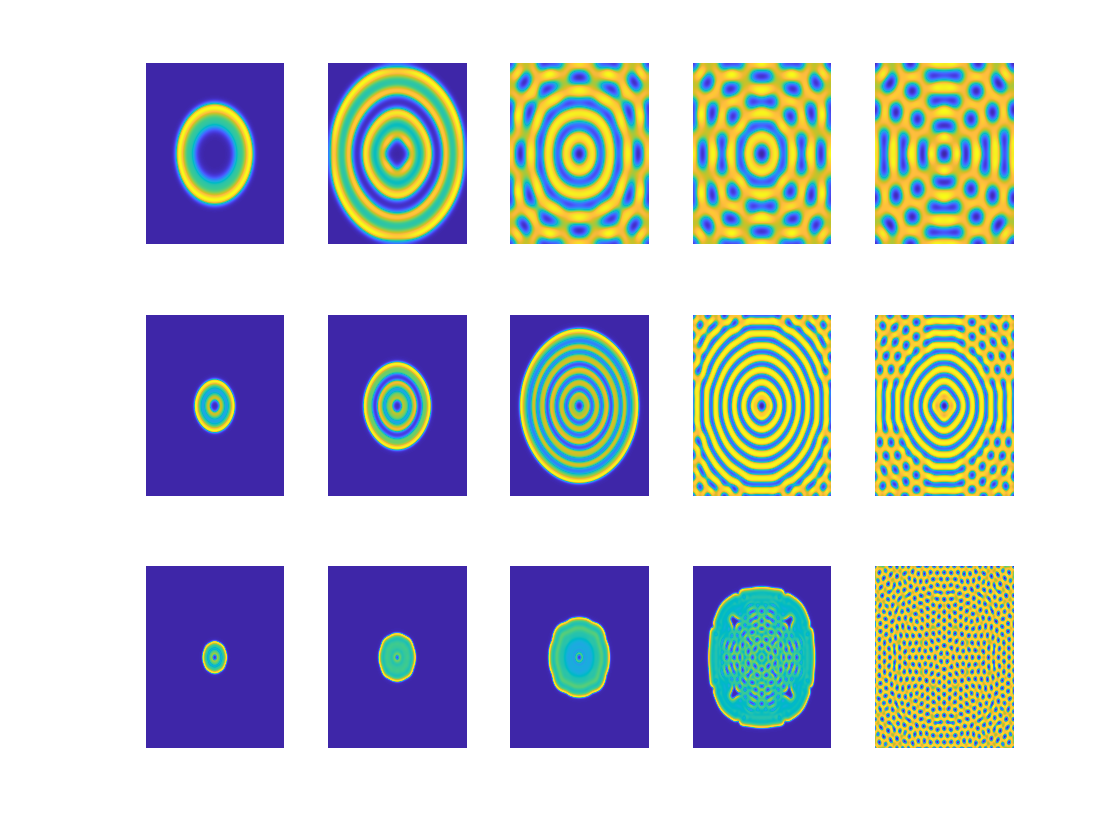}
		\caption{$ K=0.055 $, $ \alpha=2, 1.7, 1.5 $ (Up-Down), $ T=500, 1000, 2000,4000,8000 $ (Left-Right).}
		\label{fig:gsk55}
	\end{subfigure}
	~
	\begin{subfigure}[h]{1\textwidth}
		\includegraphics[height=2.2in,width=5.2in]{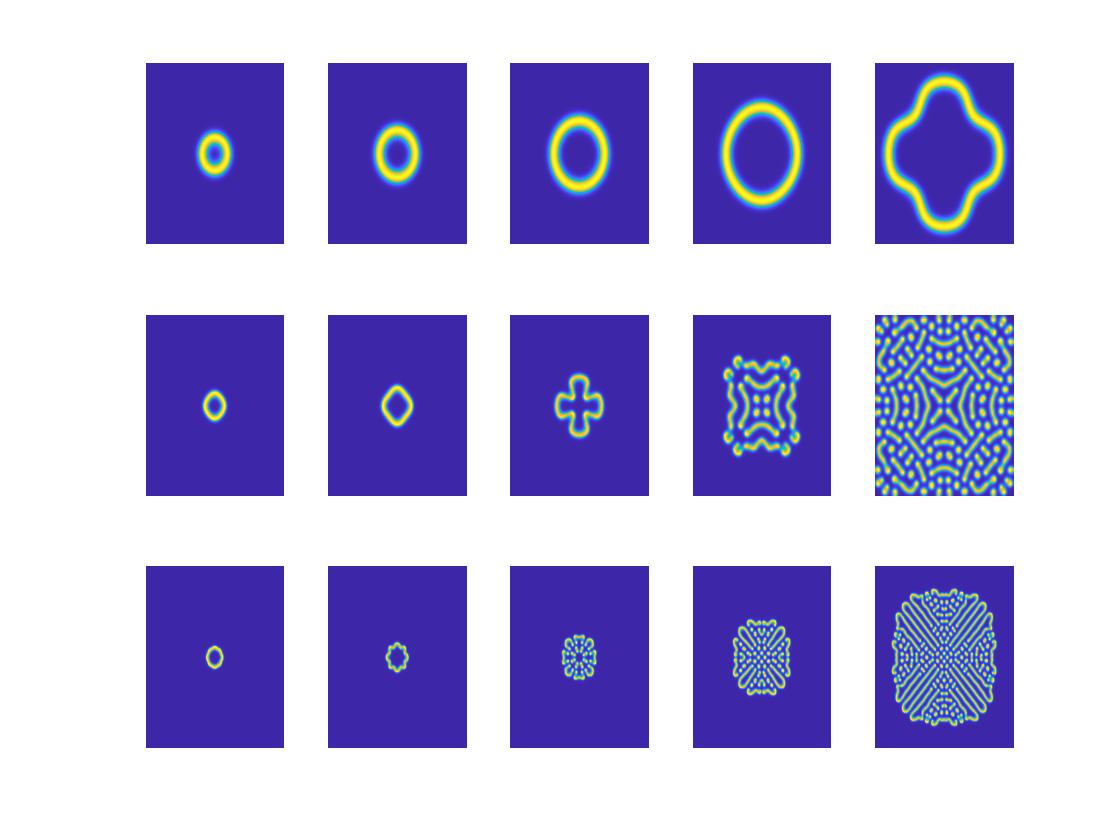}
		\caption{$ K=0.061 $, $ \alpha=2, 1.7, 1.5 $ (Up-Down), $ T=500, 1000, 2000,4000,8000 $ (Left-Right).}
		\label{fig:gsk61}
	\end{subfigure}
	~
	\begin{subfigure}[h]{1\textwidth}
		\includegraphics[height=2.2in,width=5.2in]{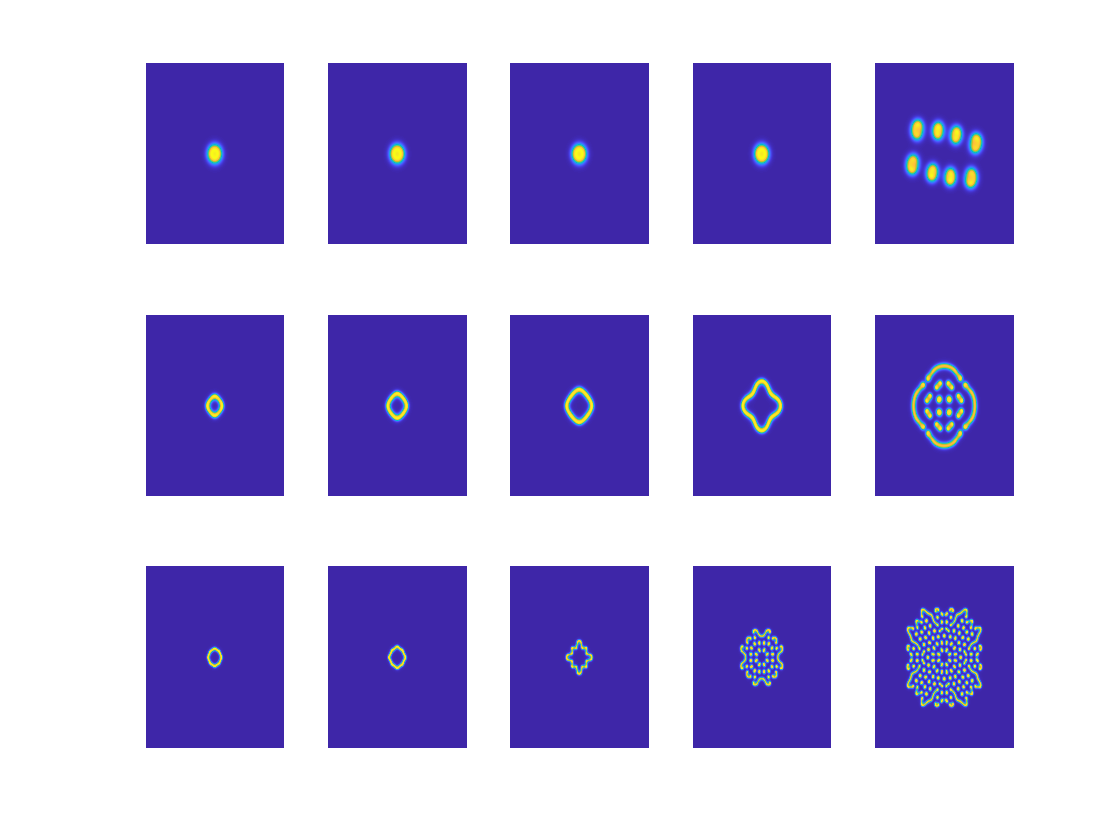}
		\caption{$ K=0.063 $, $ \alpha=2, 1.7, 1.5 $ (Up-Down), $ T=500, 1000, 2000,4000,8000 $ (Left-Right).}
		\label{fig:gsk63}
	\end{subfigure}
	\caption{Pattern formation in the Gray-Scott model (\ref{SG}) for different values of parameter $ K $ and fractional power $ \alpha $.}
	\label{fig:gs}
\end{figure}

The evolution profile of the component $ v $ and the effects of the super-diffusion for the fractional Gray-Scott model are shown in Fig. \ref{fig:gs}. From the Fig. \ref{fig:gs}, it can be seen that the speeds of pattern formulation depend on different $\alpha $ due to the speed of the diffusion is affected by fractional order.

For $  K=0.055$ (Fig. \ref{fig:gsk55}), the model with the normal diffusion ($ \alpha=2 $) produce a circular wave propagating outward. Moreover, the reduction of the fractional order $ \alpha=1.7 $ affects the size of patterns with smaller spots and generates a decrease in the velocity of the propagation of the initial perturbation. For smaller values of the fractional power ($ \alpha=1.5$), we observe a new process of pattern formation. The process of nucleation of the structure propagates outward until the entire domain reaches the final steady state. For $  K=0.061$ (Fig. \ref{fig:gsk61}), the model presents filaments and produces a wavefront propagation by curvature. When the order $ \alpha $  of the fractional Laplacian operator decreases, the curvature driven mechanisms are increased and filaments become thinner. For $  K=0.063$ (Fig. \ref{fig:gsk63}), the model exhibits patterns of the mitosis when $ \alpha=2 $. However, when the fractional order of the model is decreased, the replication pattern is completely changed as shown in Fig. \ref{fig:gsk63}.

\subsubsection{\textbf{3D Schnakenberg model}}
We consider the system of 3D space-fractional Schnakenberg model \cite{holmes2014efficient}:
\begin{align}\label{SB}
\begin{split}
\dfrac{\partial  u }{\partial  t} &= - \kappa_u (-\Delta)^{\alpha/2} u +\gamma (a-u + u^2 v),\\
\dfrac{\partial  v }{\partial  t} &= - \kappa_v (-\Delta)^{\alpha/2} v +\gamma (b-u^2v).
\end{split}
\end{align}
Subject to periodic boundary conditions and initial condition:
\begin{align*}
\begin{split}
 u(x,y,z,0) &= 1 - e^{-10 \left( (x-l/2)^2 + (y-l/2)^2 +(z-l/2)^2\right)},\\
 v(x,y,z,0) &=  e^{-10 \left( (x-l/2)^2 + 2 (y-l/2)^2 +(z-l/2)^2\right)},
 \end{split}
\end{align*}
where $ u $ and $ v $ are the concentration of the chemical products, $ \kappa_u $, $ \kappa_v $ are the diffusion coefficients, $ \gamma $ is the scale factor proportional to the length of the domain, $ a, b $ are the positive constants, and $l $ is the length of domain.
The parameters are selected as follows: $ \kappa_u=1 $, $ \kappa_v=10 $, $ \gamma=1 $, $ a=0.1 $, $ b=0.9 $. We simulate the model on the domain $ \Omega=(0,l)^3 $ with $ N=32 $ and $ \tau=1 $. 

\begin{figure}[H]
	\centering
	\includegraphics[height=7in,width=6in]{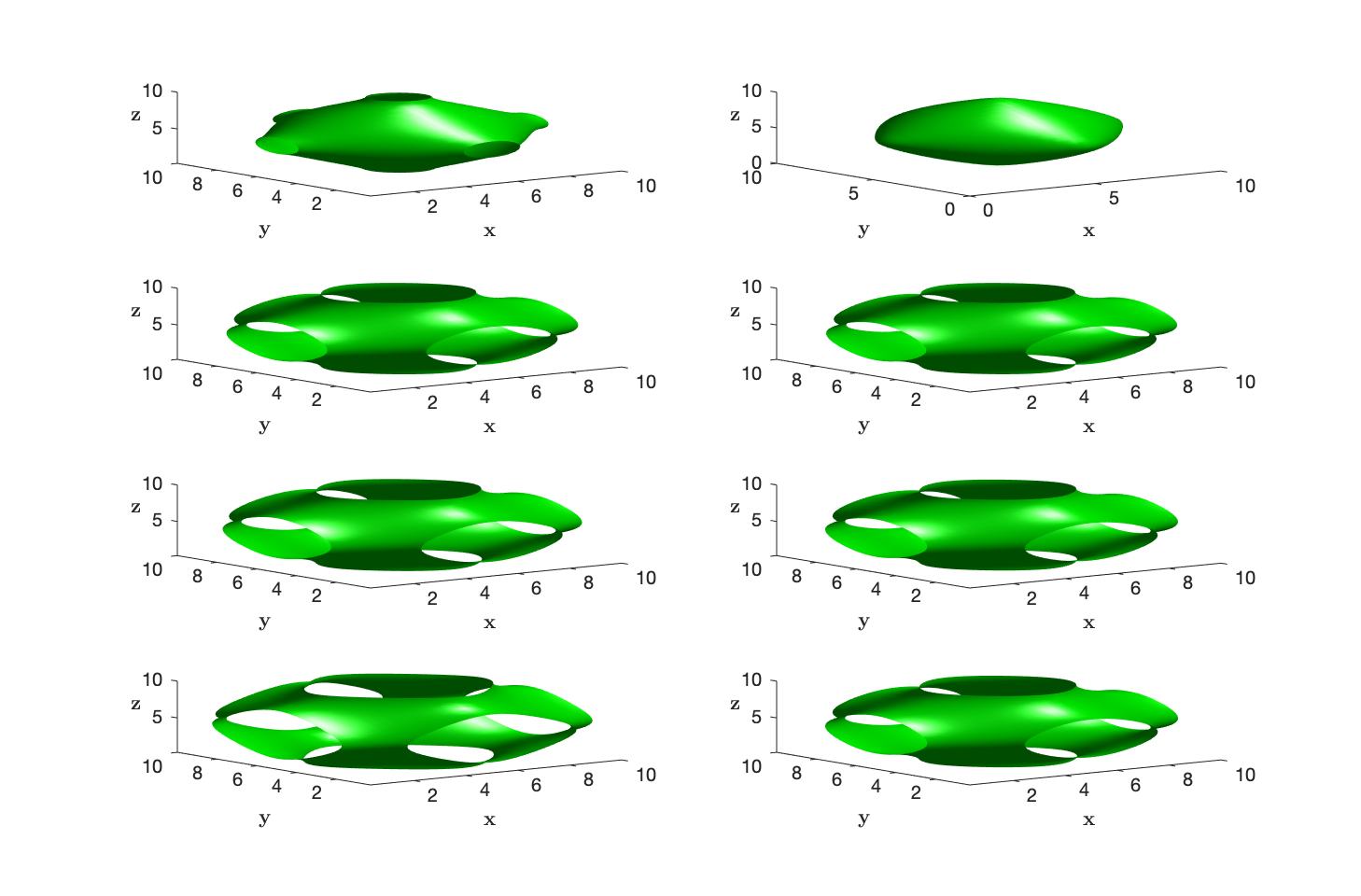}
	\caption{Snapshots of the evolution of substrate $ v $ for 3D Schnakenberg model \eqref{SB} with $ l=10 $ and $ \alpha=2,1.5 $ (Left-Right) at $ T=5, 20, 100, 200 $ (Up-Down). }
	\label{fig:sbl10}
\end{figure}

\begin{figure}[H]
	\centering
	\includegraphics[height=7in,width=6in]{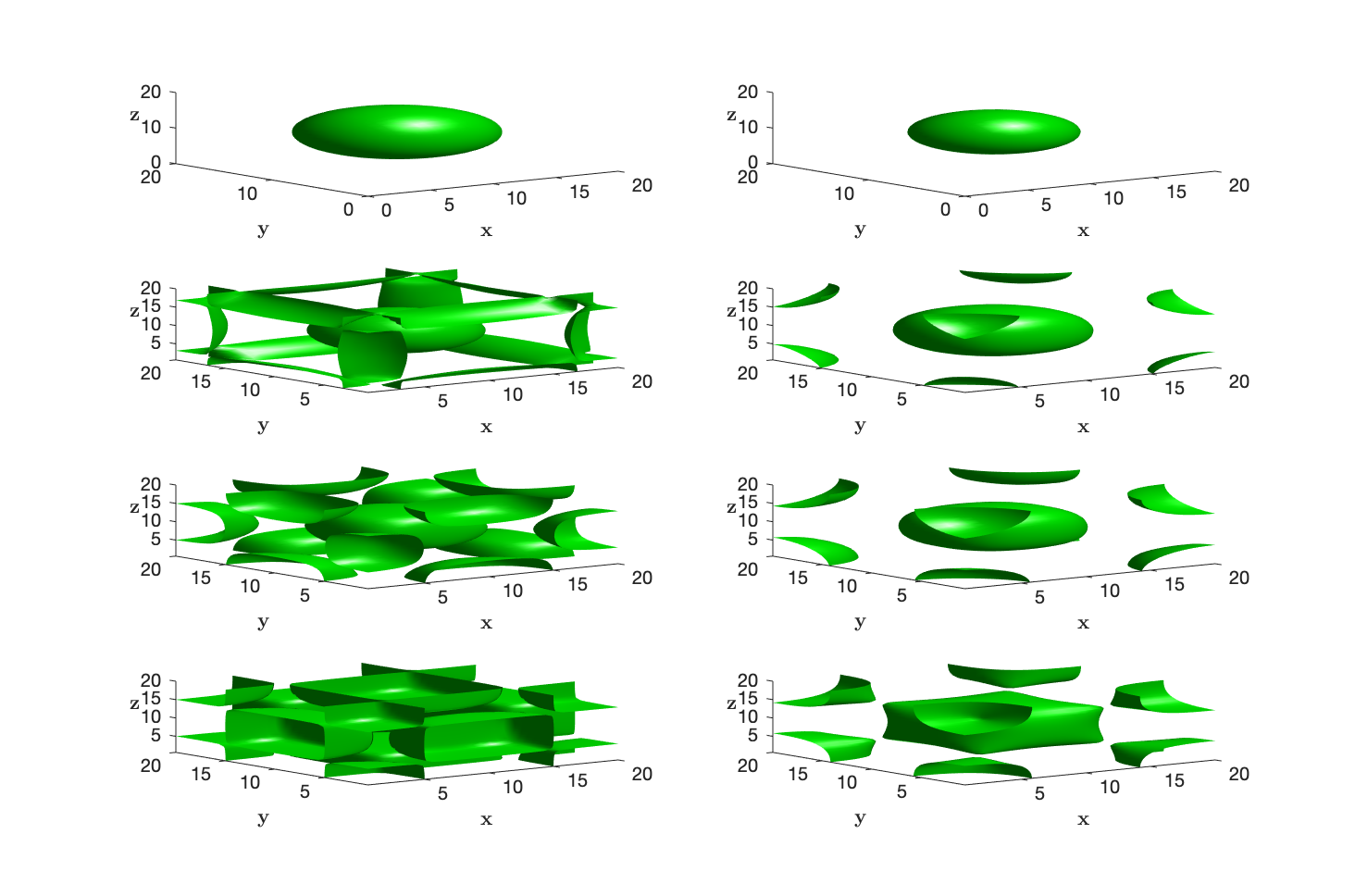}
	\caption{Snapshots of the evolution of substrate $ v $ for 3D Schnakenberg model \eqref{SB} with $ l=20 $ and $ \alpha=2,1.5 $ (Left-Right) at $ T=5, 20, 100, 200 $ (Up-Down).}
	\label{fig:sbl20}
\end{figure}

\begin{figure}[H]
	\centering
	\includegraphics[height=7in,width=6in]{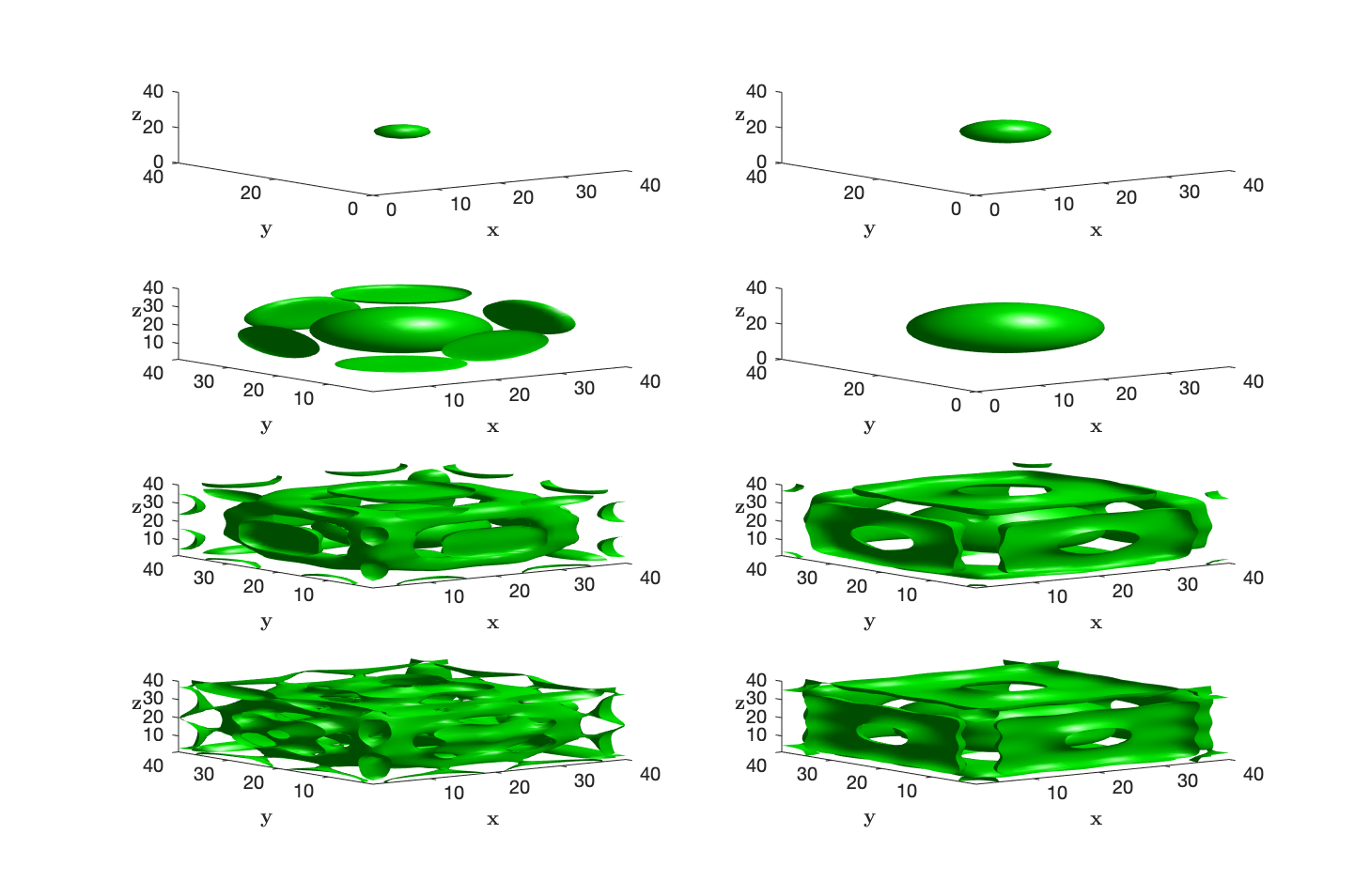}
	\caption{Snapshots of the evolution of substrate $ v $ for 3D Schnakenberg model \eqref{SB} with $ l=40 $ and $ \alpha=2,1.5 $ (Left-Right) at $ T=5, 20, 100, 200 $ (Up-Down).}
	\label{fig:sbl40}
\end{figure}

Numerical simulation of the Schnakenberg model in three spatial dimensions represents an even more amazing scenario where more exotic and chaotic patterns may arise which illustrate the emergence of pattern formation dependent on the values of $ T $, the length of the domain and the fractional power $ \alpha $, as can be observed in Figs. \ref{fig:sbl10}, \ref{fig:sbl20} and \ref{fig:sbl40}.

\section{Conclusion}
\label{S6}
%In this article, we have proposed a high-order time stepping scheme in combining with 
In this manuscript, a fast, accurate, and stable high-order numerical method is proposed for solving space-fractional reaction-diffusion equations with subject to periodic, Dirichlet or Neumann boundary conditions. The proposed method is based on a fourth-order exponential time differencing Runge-Kutta approximations for time integration and a fourth-order compact finite difference approximation in combination with matrix transfer technique for spatial approximation, that allow use of FFT-based fast calculations. An empirical convergence analysis and calculation of local truncation error exhibited the fourth-order accuracy of the proposed method. The performance (in terms of accuracy and efficiency) and reliability of the method has been investigated by testing it on various numerical examples including 2D fractional Fitzhugh-Nagumo, Gierer-Meinhardt, Gray-Scott and 3D fractional Schnakenberg models. The numerical results exhibited that the proposed method is fast, accurate and applicable for simulating the multi-dimensional space-fractional reaction-diffusion equations. In future, the method will be applied to multi-dimensional space-fractional reaction-diffusion-advection equations involving non smooth data.

\bibliographystyle{unsrt}  


%\bibliographystyle{plain}
%\bibliography{references}

\end{document}